\documentclass{elsarticle}
\usepackage{graphicx}
\usepackage{amsmath}
\usepackage{amsfonts}
\usepackage{amssymb}
\usepackage[nolist]{acronym}

\usepackage{algorithm}
\usepackage{algorithmic}

\usepackage{bm}
\usepackage{bbm}
\usepackage{url}
\usepackage{graphicx}
\usepackage{tikz}
\usepackage{pgfplots}
\usepackage{amsmath}
\usepackage{amsfonts}
\usepackage{amssymb}
\usepackage{subcaption}
\usetikzlibrary{arrows.meta, decorations.markings, intersections}
\usepgfplotslibrary{groupplots}

\begin{acronym}[mylongestacronym] 
\acro{av}[AV]{artificial viscosity}
\acro{blas}[BLAS]{Basic Linear Algebra Subprograms}
\acro{bosss}[BoSSS]{Bounded Support Spectral Solver}
\acro{cfd}[CFD]{Computational Fluid Dynamics}
\acro{cfl}[CFL]{Courant-Friedrichs-Levy}
\acro{dg}[DG]{Discontinuous Galerkin Method}
\acro{dns}[DNS]{Direct Numerical Simulation}
\acroplural{dns}[DNS]{Direct Numerical Simulations}
\acro{dof}[DOF]{degrees of freedom}
\acro{xdg}[XDG]{Extended Discontinuous Galerkin}
\acro{eoc}[EOC]{experimental order of convergence}
\acro{eos}[EOS]{equation of state}
\acrodefplural{eos}[EOS]{equations of state}
\acro{fc}[FC]{Finite Cell}
\acro{fd}[FD]{Finite Difference}
\acro{fe}[FE]{Finite Element}
\acro{fv}[FV]{Finite Volume}
\acro{gls}[GLS]{Generalized Level Set}
\acro{hllc}[HLLC]{Harten Lax van Leer Contact}
\acro{hmf}[HMF]{hierarchical moment-fitting}
\acro{hpc}[HPC]{High Performance Computing}
\acro{ibm}[IBM]{Immersed Boundary Method}
\acro{lapack}[LAPACK]{Linear Algebra Package}
\acro{lia}[LIA]{linear interaction analysis}
\acro{lts}[LTS]{local time stepping}
\acro{ode}[ODE]{ordinary differential equation}
\acro{pde}[PDE]{partial differential equation}
\acrodefplural{pde}[PDEs]{partial differential equations}
\acro{simple}[SIMPLE]{Semi-Implicit Method for Pressure-Linked Equations}
\acro{tvd}[TVD]{total variation deminishing}
\acro{xfem}[XFEM]{eXtended Finite Element Method}
\acro{sti}[STI]{shock-turbulence interaction}
\acro{1d}[1D]{one-dimensional}
\acro{2d}[2D]{two-dimensional}
\acro{3d}[3D]{three-dimensional}
\end{acronym}
\usepackage{amsfonts} 

\newcommand{\jump}[1]{\left[\!\left[{#1}\right]\!\right]}

\newcommand{\nbrc}[1]{[#1]}
\newcommand{\suchthat}{\mathrel{}\middle|\mathrel{}}

\newcommand{\Mach}{\mathrm{Ma}}

\newcommand*{\termn}{n_{\text{term}}}
\newcommand{\inn}[1]{#1^{\text{in}}}
\newcommand{\out}[1]{#1^{\text{out}}}
\newcommand{\pder}[2]{\ensuremath{\frac{\partial #1}{\partial #2}}}

\newcommand{\Acal}{\ensuremath{\mathcal{A}}}

\newcommand{\Kcal}{\ensuremath{\mathcal{K}}}

\newcommand{\Vcal}{\ensuremath{\mathcal{V}}}

\newcommand{\Nbb}{\ensuremath{\mathbb{N} }}

\newcommand{\Rbb}{\ensuremath{\mathbb{R} }}

\newcommand\Abm{{\ensuremath{\bm{A}}}}
\newcommand\Bbm{{\ensuremath{\bm{B}}}}

\newcommand\Hbm{{\ensuremath{\bm{H}}}}
\newcommand\Ibm{{\ensuremath{\bm{I}}}}
\newcommand\Jbm{{\ensuremath{\bm{J}}}}

\newcommand\Rbm{{\ensuremath{\bm{R}}}}

\newcommand\cbm{{\ensuremath{\bm{c}}}}

\newcommand\fbm{{\ensuremath{\bm{f}}}}

\newcommand\nbm{{\ensuremath{\bm{n}}}}

\newcommand\rbm{{\ensuremath{\bm{r}}}}

\newcommand\ubm{{\ensuremath{\bm{u}}}}
\newcommand\vbm{{\ensuremath{\bm{v}}}}

\newcommand\xbm{{\ensuremath{\bm{x}}}}

\newcommand\zbm{{\ensuremath{\bm{z}}}}

\newcommand\lambdabold{{\ensuremath{\boldsymbol{\lambda}}}}

\newcommand\phibold{{\ensuremath{\boldsymbol{\phi}}}}
\newcommand\varphibold{{\ensuremath{\boldsymbol{\varphi}}}}

\newcommand\Phibold{{\ensuremath{\boldsymbol{\Phi}}}}

\title{An extended discontinuous Galerkin shock tracking method}
\author[1,2]{Jakob Vandergrift\corref{cor1}}
\ead{vandergrift@fdy.tu-darmstadt.de}

\author[1,2]{Florian Kummer}
\ead{kummer@fdy.tu-darmstadt.de}

\cortext[cor1]{Corresponding author}
\affiliation[1]{organization={TU Darmstadt, Chair of Fluid Dynamics},\\ 
	addressline={Otto-Berndt-Str. 2 },
	postcode={64287}, 
	city={Darmstadt}, 
	country={Germany}}
\affiliation[2]{organization={TU Darmstadt, Graduate School of Computational Engineering},\\
	addressline={Dolivostraße 15},
	postcode={64293}, 
	city={Darmstadt}, 
	country={Germany}}

\begin{document}
	\begin{abstract}
		In this paper, we introduce a novel high-order shock tracking method and provide a proof of concept. Our method leverages concepts from implicit shock tracking and extended discontinuous Galerkin methods, primarily designed for solving partial differential equations featuring discontinuities. To address this challenge, we solve a constrained optimization problem aiming at accurately fitting  the zero iso-contour of a level set function to the discontinuities. Additionally, we discuss various robustness measures inspired by both numerical experiments and existing literature. Finally, we showcase the capabilities of our method through a series of two-dimensional problems, progressively increasing in complexity.
	\end{abstract}
	\maketitle

\section{Introduction}
In the field of Computational Fluid Dynamics (CFD), high-order methods, such as the Discontinuous Galerkin method (DG), are under current research, as they provide several benefits for accurately simulating fluid flow at a high level of fidelity including high accuracy per \ac{dof}, parallel scalability and low dissipation. However, when simulating transonic compressible flows, a numerical method must also be able to handle discontinuities in the flow field arising from shock waves. 
Shock waves are a challenging phenomenon in transonic flows arising whenever the velocities exceed the local speed of sound ($\Mach>1$) \cite{andersonModernCompressibleFlow2003}. The transition from super- to subsonic flow forms a thin layer in the flow field, where flow properties change drastically. 
These thin layers can be modelled as discontinuity surfaces along which scalar quantities jump. In their vicinity and without any stabilization measures, high-order methods suffer from non-linear instabilities, leading to decreases in accuracy and ultimately failure of the simulation. 

As a result, various stabilization techniques have been devised, with many relying on \textit{shock-capturing} strategies. In such approaches, the numerical discretization is tailored to handle discontinuities, especially in cells where oscillations in the solution are detected. These techniques encompass the use of limiters to attain total variation diminishing schemes \cite{cockburnRungeKuttaDiscontinuous2001}, high-order reconstruction methods like the weighted essentially non-oscillatory (WENO) scheme \cite{hartenUniformlyHighOrder1997,shuHighOrderWeighted2009}, and the introduction of artificial viscosity \cite{perssonSubCellShockCapturing2006,barterShockCapturingPDEbased2010,chingShockCapturingDiscontinuous2019}. While these approaches have demonstrated their effectiveness in practice, shock-capturing strategies do come with certain drawbacks. They tend to dampen low-order perturbations, including acoustic waves. Moreover, in close proximity to shock fronts, their accuracy diminishes, often degrading to first-order accuracy \cite{wangHighorderCFDMethods2013}.

A less conventional approach to handling shock waves involves aligning the edges of computational mesh cells with the shock front, a technique commonly referred to as \textit{shock fitting/tracking}. By leveraging inter-element jumps in numerical solutions, these methods allow for the exact representation of discontinuities without the need for additional stabilization techniques. However, this approach comes with its own set of challenges. One of the primary challenges is that the precise position of the shock wave is often unknown, and it can change or evolve over time, making meshing a complex and dynamic problem. Additionally, complex patterns can emerge as a result of interactions between different shock waves or reflections. Early approaches, such as those mentioned in \cite{trepanierConservativeShockFitting1996} and \cite{bainesMultidimensionalLeastSquares2002}, primarily focused on low-order schemes where shock capturing had a more significant relative advantage. Consequently, shock fitting has played a relatively minor role in the simulation of transonic flows until more recently. One class of the methods which are still under current research \cite{assonitisNewShockfittingTechnique2022} employs specialized strategies to explicitly identify the shock front using the Rankine-Hugoniot conditions and are surveyed in \cite{salasShockFittingPrimer2009,morettiThirtysixYearsShock2002}. While these approaches have shown promising results, they may not be easily adaptable to problems where the shock's position and topology are unknown or subject to change.

In recent years, new \textit{implicit shock fitting/tracking} approaches have emerged within the context of DG methods and shown promise in addressing the challenges associated with shock waves in supersonic flows. These methods treat mesh coordinates as additional variables within the discretized conservation law and compute their numerical solutions, implicitly aligning the mesh edges with discontinuities. In the work of Corrigan et al.  \cite{corriganMovingDiscontinuousGalerkin2019,corriganConvergenceMovingDiscontinuous2019a}, the Moving Discontinuous Galerkin method with Interface Condition Enforcement (MDG-ICE) is proposed where the corresponding weak form of a system of stationary conservation laws is augmented by a term enforcing the Rankine-Hugoniot conditions. The resulting non-linear system is then solved by a Levenberg-Marquardt method. Subsequently, a least squares formulation of the same method was presented \cite{kercherLeastsquaresFormulationMoving2021} and extended to handle viscous flows \cite{kercherMovingDiscontinuousGalerkin2021a}. 
Another significant development is the High-Order Implicit Shock Tracking (HOIST) method by Zahr et al. \cite{zahrRAdaptiveHighOrderDiscontinuous2020a}. HOIST achieves shock fitting of the mesh by solving a constrained optimization problem, a strategy that is also employed in the present work. Furthermore, this method has been extended to address reactive flows \cite{zahrHighOrderResolutionMultidimensional2021}, time-dependent problems \cite{shiImplicitShockTracking2022} and additional robustness was investigated along with steady 3D cases \cite{huangRobustHighorderImplicit2022a}.

The main goal of our work is to apply the eXtended Discontinuous Galerkin (XDG) method, also known as the unfitted DG (UDG) method, in a shock tracking framework. In the XDG method, level set functions are employed to represent immersed boundaries and surfaces of discontinuities implicitly by their zero iso-contours. The latter effectively intersect the fixed computational background grid, splitting/cutting its cells into sub-cells. Within these cut-cells, the traditional DG basis functions are replaced with XDG basis functions, enabling the accurate representation of discontinuities at the implicitly defined interfaces. The XDG method has previously been employed to represent boundaries of the computational domain alongside shock capturing strategies for simulating compressible flow \cite{mullerHighorderDiscontinuousGalerkin2017,geisenhoferDiscontinuousGalerkinImmersed2019}, while mainly demonstrating its effectiveness in other areas of CFD that involve discontinuous solutions, including multi-phase flows \cite{kummerExtendedDiscontinuousGalerkin2017a,smudaMarchingLevelSet2022}. In the context of explicit shock fitting, a reconstruction algorithm has been proposed \cite{geisenhofermarkusShockCapturingHighOrderShockFitting2021} with the goal of recovering the aligned level set from prior shock capturing simulations. However, the reconstructed level set was observed to be only sub-cell accurate and did not stabilize the underlying DG method properly. 

In this work, we introduce a novel implicit XDG shock tracking method. The method aims at computing shock-aligned XDG solutions to supersonic flow problems and does not require stabilization through shock capturing. To align the level set accurately with the shock, we treat its coefficients as additional variables in the discretized conservation laws and formulate an optimization problem adapted from Zahr et. al \cite{zahrImplicitShockTracking2020}. When compared to mesh-based shock tracking, this approach offers some advantages, including the ability to keep the grid fixed, thereby avoiding computational challenges associated with mesh operations, such as handling ill-shaped cells. Mainly, the method holds promise for problems involving moving shocks, as it can leverage existing techniques for evolving interfaces in time \cite{kummerTimeIntegrationExtended2018}. Also, traditional shock tracking methods might necessitate re-meshing when dealing with moving shocks, a step that can be circumvented using our approach.

This work is structured as follows: First, DG and XDG approximation spaces are introduced in Section 2. Secondly, in Section 3, these spaces are used to discretize the generic system of conservation laws and formulate a constrained optimization problem, aiming at level set shock-alignment. Consequently, a Sequential Quadratic Programming (SQP) method for solving the problem is developed in Section 4. In Section 5, an explicit parametrization of the level set function, needed to enhance convergence of the shock tracking method, is discussed. Section 6 gives details on robustness measures incorporated into the shock tracking framework and is followed by Section 7, where numerical results are discussed. Finally, the paper is concluded and an outlook for future work is given in Section 8.
\subsection{BoSSS - open source software}
The method described in this work was developed and implemented as part of the open-source software framework BoSSS \cite{bosss}, short for Bounded Support Spectral Solver. BoSSS offers a flexible environment tailored to the development, evaluation, and application of numerical discretization strategies for partial differential equations (PDEs), with the XDG method as its main research subject. The BoSSS code has found successful applications in a range of challenging scenarios. These applications include simulations of incompressible flows with active particles \cite{deussenHighorderSimulationScheme2021}, compressible inviscid and viscous flows with moving immersed boundaries \cite{mullerHighorderDiscontinuousGalerkin2017}, extension for shock-capturing in high Mach number flows \cite{geisenhoferDiscontinuousGalerkinImmersed2019}, diffusion flames in low Mach number flows \cite{gutierrez-jorqueraFullyCoupledHigh2022}, incompressible multiphase flows \cite{Kummer2016}, viscoelastic fluid flows \cite{kikkerFullyCoupledHigh2021} and helical flows \cite{dierkesHighOrderDiscontinuousGalerkin2022}.

All XDG shock tracking test cases discussed in this work (Section 7) can be reproduced by downloading the publicly available repository\footnote{
	\url{https://github.com/FDYdarmstadt/BoSSS}} and starting the respective jupyter notebooks from the directory \url{./examples/ShockFitting}.
\section{Preliminaries}
In this section, we define the approximation space in which we aim to numerically solve hyperbolic PDEs with discontinuities in the solution, e.g. the Euler equations in the vicinity of shocks. We consider a general system of $m \in \{1,\hdots,4\}$ inviscid conservation laws defined on a $D$-dimensional (in this work $D=2$) physical domain $\Omega \subset \mathbb{R}^D$
\begin{eqnarray}
	\label{equation:conservationlaw}
	\nabla \cdot \fbm(U) =0 \quad \text{in } \Omega,
\end{eqnarray}
subject to appropriate boundary conditions.
Here, $U:\Omega \rightarrow \mathbb{R}^m$ denotes the solution of the system, $\nabla := (\partial_{x_1},\dots,\partial_{x_D})$ the gradient operator in the physical domain and $\fbm:\mathbb{R}^m \rightarrow \mathbb{R}^{m\times D}$ denotes the physical flux. 

\paragraph{Discontinuous Galerkin space} The domain $\Omega$ is discretized using a Cartesian background grid $\mathcal{K}_h$ defined as 
\begin{eqnarray}
	\mathcal{K}_h = \{ K_1 ,K_2 , \ldots, K_J \},
\end{eqnarray}
covering the whole domain with non-overlapping cells, e.g. $\bigsqcup_{i=1}^J K_i = \Omega$.
Furthermore, we denote by $\nbm_\Gamma$ the normal field on the set of all edges $\Gamma = \bigcup_{j\in J} \partial K_j$. One can then define the vector-valued DG space
\begin{eqnarray}
	\mathcal{V}^m_P := \lbrack \mathbb{P}_P(\mathcal{K}_h)\rbrack ^m =  \left\{V \in \lbrack L^2(\Omega) \rbrack ^m \suchthat V\vert_{K} \in \lbrack \mathbb{P}_P(K)\rbrack ^ m, \forall K \in \mathcal{K}_h \right\} ,
\end{eqnarray}
where $\mathbb{P}_P(K)$ denotes the space of polynomials of degree $P \in \mathbb{N}$ defined on the cell $K\in \mathcal{K}_h$. For a field $\psi \in C^0(\Omega \setminus \Gamma) \subset  \mathcal{V}_P :=\mathcal{V}^1_P$, we denote by  
\begin{eqnarray}
	\inn{\psi}:= \lim_{\epsilon \rightarrow 0^+} \psi(\xbm-\epsilon \nbm_\Gamma), \quad \out{\psi}:= \lim_{\epsilon \rightarrow 0^-} \psi(\xbm-\epsilon \nbm_\Gamma), \quad \forall \xbm \in \Gamma
\end{eqnarray}
the inner and outer traces and the jump operator by
\begin{eqnarray}
	\jump{\psi}:= \begin{cases} \inn{\psi}-\out{\psi} \quad &\text{on }\Gamma_{\text{int}} \\
		\inn{\psi} \quad &\text{on }\partial\Omega 
	\end{cases}.
\end{eqnarray}
 Furthermore, we introduce the notation $\nbrc{N}:=\{1,2,\hdots,N\}$ for natural numbers $N\in \Nbb$. Then, for every cell \( K_j \) we denote by \( \{\phi_{j,n}\}_{n \in \nbrc{N_P} } \) a modal basis of the space $\mathbb{P}_P(K_j)$ and introduce the polynomial basis vector $\Phibold(P)=\Phibold $ by
\begin{eqnarray}
	\label{eqn:disc:dg:basisvec}
	\phi_{j} = (\phi_{j,n})_{n \in \nbrc{N_P}} \in (\mathcal{V}_P)^{N_P} ~,~\Phibold:= (\phi_{j})_{j \in \nbrc{J}}, \quad \Phibold^m=  
	(\underbrace{\Phibold
		~\hdots~
		\Phibold }_{m \text{ times}})^T.
\end{eqnarray}
Then, each element \( V\in \mathcal{V}^m_P \) residing in the DG space can be described by an unique coordinate vector by
\begin{eqnarray}
		\vbm
		\in \mathbb{R}^{m\cdot J \cdot N_P} \text{ s.t. } V=\vbm \cdot \Phibold.
	\end{eqnarray}
	\paragraph{Extended discontinuous Galerkin space} Next, the DG method will be extended using two level sets $\varphi_{b}\in \mathcal{V}_{P_b},~\varphi_{s} \in \mathcal{V}_{P_s}$ with respective polynomial degrees $P_b, P_s \geq 1$. These level sets implicitly influence the discretization by their zero-sets, either representing discontinuities in the solution ($\{\varphi_s = 0\}$) or boundaries ($\{\varphi_b=0\}$) immersed into the background grid $\mathcal{K}_h$. They implicitly decompose the domain $\Omega$ into (level set dependent) sub-domains (SD) $\mathfrak{A},\mathfrak{B},\mathfrak{C},\mathfrak{D}$  which are defined by (see Figure \ref{fig:subdomains} for illustration)
\begin{align} 
	\label{eq:subdomains}
	\mathfrak{A} = \mathfrak{A}(\varphi_s,\varphi_b) &:=  \left\{\xbm\in \Omega \suchthat \varphi_s (\xbm) <0 ~\cap ~\varphi_b (\xbm) <0\right\} , \nonumber \\
	\mathfrak{I}_s = \mathfrak{I}(\varphi_s) &:=\left\{\xbm\in \Omega \suchthat \varphi_s (\xbm) =0  \right\} , \nonumber\\
	\mathfrak{B} = \mathfrak{B}(\varphi_s,\varphi_b) &:=\left\{\xbm\in \Omega \suchthat \varphi_s (\xbm) >0 ~\cap ~\varphi_b (\xbm) <0 \right\}, \nonumber\\
	\mathfrak{C} = \mathfrak{C}(\varphi_s,\varphi_b) &:=\left\{\xbm\in \Omega \suchthat \varphi_s (\xbm) <0 ~\cap ~\varphi_b (\xbm) >0 \right\}, \\
	\mathfrak{I}_b = \mathfrak{I}(\varphi_b) &:=\left\{\xbm\in \Omega \suchthat \varphi_b (\xbm) =0 \right\}, \nonumber\\
	\mathfrak{D}=\mathfrak{D}(\varphi_s,\varphi_b)  &:=\left\{\xbm\in \Omega \suchthat \varphi_s (\xbm) >0 ~\cap ~\varphi_b (\xbm) >0 \right\}, \nonumber\\
	\text{SD}=\text{SD}(\varphi_s,\varphi_b)&:=\left\{\mathfrak{A},\mathfrak{B} ,\mathfrak{C},\mathfrak{D}\right\}. \nonumber
\end{align}
	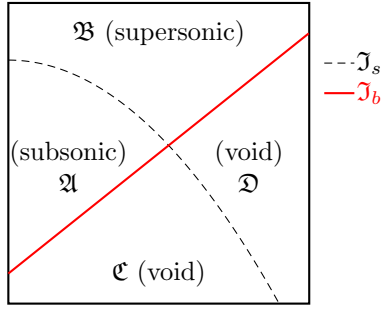
\begin{figure}
	\centering
	\begin{tikzpicture}[>=Stealth, scale=4]
		\draw[thick] (0,0) rectangle (1,1);
		
		\draw[black, densely dashed, domain=0:0.9, samples=100, name path=curve1] plot (\x, {0.81 - \x*\x});
		
		\draw[red, thick, domain=0:1, samples=100, name path=curve2] plot (\x, {0.1 + 0.8*\x});
		
		\node at (0.2,0.4) {$\mathfrak{A}$};
		\node at (0.2,0.5) {(subsonic)};
		\node at (0.8,0.4) {$\mathfrak{D}$};
		\node at (0.8,0.5) {(void)};
		\node at (0.5,0.9) {$\mathfrak{B}$  (supersonic)};
		\node at (0.5,0.1) {$\mathfrak{C}$  (void)};
		
		\node[black] at (1.2,0.8) {$\mathfrak{I}_s$};
		\draw[black, densely dashed] (1.05,0.8) -- (1.15,0.8);
		\node[red] at (1.2,0.7) {$\mathfrak{I}_b$};
		\draw[red, thick] (1.05,0.7) -- (1.15,0.7);
		
	\end{tikzpicture}
	\caption{Exemplary cell given two interfaces $\mathfrak{I}_s,\mathfrak{I}_b$ (defined by $\varphi_s$ and $\varphi_b$). The resulting sub-domains $\mathfrak{A}$,$\mathfrak{B}$,$\mathfrak{C}$,$\mathfrak{D}$ are depicted as defined in \eqref{eq:subdomains}. Also, exemplary flow regimes which are used for supersonic flow configurations examined in this work (\textit{supersonic,~subsonic,~void}) are assigned to each sub-domain.}
	\label{fig:subdomains}
\end{figure}
To obtain the computational grid, which is used for the XDG discretization, each cell $K_j$, intersected by one of the interfaces ($\mathfrak{I}_b,	\mathfrak{I}_s$) is split/cut into cut-cells $K_{j,\mathfrak{s}}:=  K_j \cap \mathfrak{s}$. The resulting level set dependent cut-cell grid $\mathcal{K}^X_h$ is then defined as
\begin{eqnarray}
	\mathcal{K}^X_h=\mathcal{K}^X_h(\varphi_s,\varphi_b) := \left\{ K_{j,\mathfrak{s}} :=  K_j \cap \mathfrak{s} \suchthat j\in\nbrc{J}, \forall\mathfrak{s}\in \text{SD} \right\},
\end{eqnarray}
and using the above, the vector-valued XDG space can be defined as
\begin{eqnarray}
	\mathcal{V}^{X,m}_P := \lbrack \mathbb{P}^X_P(\mathcal{K}_h)\rbrack ^m = \left\{ V\in \lbrack L^2(\Omega) \rbrack ^m \suchthat V\vert_{K \cap \mathfrak{s}} \in \lbrack \mathbb{P}_P(K \cap \mathfrak{s})\rbrack ^ m, \forall K \in \mathcal{K}_h,\forall\mathfrak{s} \in \text{SD} \right\}.
\end{eqnarray}
We aim to approximate the solution $U$ of \eqref{equation:conservationlaw} by functions $U_h$ residing in the XDG space $\mathcal{V}^{X,m}_P$.
To express the solution in terms of coefficients, we first define a respective basis for each sub-domain $ \mathfrak{s} \in \text{SD}$ by multiplying every element of the DG basis by an indicator function $\mathbbm{1}_{\mathfrak{s}}$ supported in the corresponding sub-domain $\mathfrak{s}$. In this sense, the XDG basis vector $\Phibold_{\mathfrak{s}}(P)$ for one sub-domain is obtained using the DG basis vector $\Phibold(P)$ (see \eqref{eqn:disc:dg:basisvec}) and element-wise multiplication with $\mathbbm{1}_{\mathfrak{s}}$, i.e.
\begin{eqnarray}
	\Phibold_{\mathfrak{s}}=\Phibold_{\mathfrak{s}}(P)= \Phibold(P) \cdot \mathbbm{1}_{\mathfrak{s}}.
\end{eqnarray}
Note, that the resulting XDG basis functions $\phi_{j,n,\mathfrak{s}}=\phi_{j,n}\mathbbm{1}_{\mathfrak{s}}$ are only supported in the respective cut-cells $K_{j,\mathfrak{s}}$, hence allowing for additional discontinuities along the interfaces $\mathfrak{I}_b,\mathfrak{I}_s$. 
	Furthermore, we introduce a full XDG basis vector $\Phibold^m_X$ for the space $\mathcal{V}^{X,m}_P$ containing all subdomains and $m$ components by 
	\begin{eqnarray}
		\Phibold_X=
		\begin{pmatrix}
			\Phibold_\mathfrak{A} \\
			\Phibold_\mathfrak{B} \\
			\Phibold_\mathfrak{C} \\
			\Phibold_\mathfrak{D} \\
		\end{pmatrix}, \quad \text{and } \Phibold^m_X=
		(\underbrace{\Phibold_X
			~\hdots~
			\Phibold_X }_{m \text{ times}})^T,
	\end{eqnarray}
	and a coefficient vector $\vbm$ for each element \( V\in \mathcal{V}^{X,m}_P \) defined by
	\begin{eqnarray} 	
		\vbm \in \mathbb{R}^{N_u}  \text{ s.t. } V=\vbm \cdot \Phibold^m_X,
	\end{eqnarray}
 	where $N_u= 4m\cdot J \cdot N_P= \dim(\mathcal{V}^{X,m}_P)$. It is important to note, that this definition of the basis vector introduces numerous zero-valued basis functions in cells that are not intersected by both interfaces, commonly referred to as ghost cells. From a mathematical perspective, these zero functions need to be removed to establish a well-defined unique basis. 
	\section{XDG discretization and problem formulation}
 In this section, the XDG discretization of the general system of conservation laws \eqref{equation:conservationlaw} is provided. On its basis, we formulate an optimization problem which is at the heart of the XDG shock tracking method. 
 
 \paragraph{XDG discretization} To obtain an XDG discretization of \eqref{equation:conservationlaw}, the Galerkin approach is used, multiplying every component of the target system \eqref{equation:conservationlaw} by suitable XDG test functions \( v \in \mathcal{V}^{X,1}_P \) and integrating by parts to obtain the weak form
	\begin{eqnarray}
		\label{equation:weakform}
		\oint_\Gamma \fbm_i(U) \cdot \nbm_\Gamma \jump{v} ~dS  - \int_\Omega \fbm_i(U) \cdot \nabla_h v ~dV=0. 
	\end{eqnarray}
The solution $U$ is approximated in the XDG space $\mathcal{V}^{X,m}_P$ by 
	\begin{eqnarray}
		\label{eqn:U_approx}
		U \approx U_h := \ubm \cdot \Phibold_X,
	\end{eqnarray}
	and is inserted into equation \eqref{equation:weakform}. Due to the ambiguous values of $U_h$ on cell edges $\Gamma$ a numerical flux 
	\begin{eqnarray}
		\label{eqn:numflux}
		\hat{F}_i(U^+,U^-,\nbm_\Gamma) \approx \fbm_i(U) \cdot \nbm_\Gamma
	\end{eqnarray}
	must be introduced for every flux component $\fbm_i$. Inserting a specific basis function $ \phi_{j,n,\mathfrak{s}} $ and \eqref{eqn:numflux} into \eqref{equation:weakform}, we obtain the respective cut-cell local residual $ 	r^{i,j,\mathfrak{s}}_{n}: \mathcal{V}^{X,m}_P \times \mathcal{V}_{P_b} \times \mathcal{V}_{P_s} \rightarrow \mathbb{R}$ which writes
	\begin{eqnarray}
		\label{eqn:localresidual}
		r^{i,j,\mathfrak{s}}_{n}(U_h,\varphi_s,\varphi_b) := \oint_{\partial K_{j,\mathfrak{s}}} \hat{F}_i(U_h^+, U_h^-,\nbm_{\Gamma}) \jump{\phi_{j,n,\mathfrak{s}}} dS - \int_{K_{j,\mathfrak{s}}}\fbm_i(U_h) \cdot \nabla_h \phi_{j,n,\mathfrak{s}}  dV.
	\end{eqnarray}
	Collecting all of these into one array, the full residual $ \rbm_P: \mathcal{V}^{X,m}_P \times \mathcal{V}_{P_b} \times \mathcal{V}_{P_s} \rightarrow \Rbb^{N_U}$ is defined and we obtain a system of equations
	\begin{eqnarray}
		\rbm_P(U_h,\varphi_s,\varphi_b):= \begin{pmatrix}
			\vdots \\
			r^{i,j,\mathfrak{s}}_{n}(U_h,\varphi_s,\varphi_b) \\
			\vdots
		\end{pmatrix}=0.
	\end{eqnarray}
	To obtain an algebraic system representing the discretized system of conservation laws, we introduce the algebraic residual $ \rbm: \Rbb^{N_U} \times \Rbb^{N_{s}} \rightarrow \mathbb{R}^{N_U}$ as
	\begin{eqnarray}
		\label{eqn:algebraic_residual}
		\rbm(\ubm,\varphibold):=\rbm_P(\ubm \cdot \Phibold_X,l_s(\varphibold) \cdot \Phibold(P_s), \varphi_b),
	\end{eqnarray}
fixing the polynomial degree $P$. Note, that in this study, we are only interested in optimizing/fitting the shock interface $\mathfrak{I}_s=\{ \varphi_s=0\}$. We thus consider the immersed boundary $\mathfrak{I}_b$ and its associated level set $\varphi_b$ as fixed. Hence, the dependence on $\varphi_b$ is discarded.

	\paragraph{Enriched XDG discretization} We seek to formulate an optimization problem, which inherently enforces alignment of the interface $\mathfrak{I}_s$ with discontinuities. As demonstrated in \cite{corriganMovingDiscontinuousGalerkin2019}, the residual of a DG method can vanish without properly fitting the interface. Therefore, it is not employed as the objective function in the optimization problem but rather as a constraint. For the objective function we opt for an enriched residual of higher degree denoted as $P^{*}>P$ (choosing $P^{*}=P+1$ in this work). This choice introduces additional constraints and is more sensitive to non-aligned interfaces and oscillatory solutions, a fact discussed in detail in \cite{zahrImplicitShockTracking2020}. To define the enriched residual $\Rbm$ in terms of $\rbm_{P^{*}}$, consider $U_h \in  \mathcal{V}^{X,m}_P $ and a canonical injection $i: \mathcal{V}^{X,m}_P \rightarrow  \mathcal{V}^{X,m}_{P^{*}}$. Then, one can define the algebraic enriched residual $\Rbm: \mathbb{R}^{N_U} \times \mathbb{R}^{N_s} \rightarrow \mathbb{R}^{4 \cdot m \cdot N_{P^*} \cdot J}$ using the residual $\rbm_{P^{*}}$ by
	\begin{eqnarray}
		\Rbm(\ubm,\varphibold) =\rbm_{P^{*}}(i(\ubm \cdot \Phibold_X),l_s(\varphibold) \cdot \Phibold(P_s), \varphi_b).
	\end{eqnarray}
Comparing to $\rbm(\ubm,\varphibold)$, the enriched residual $\Rbm(\ubm,\varphibold)$ features more function components as additional basis functions corresponding to $\mathcal{V}^{X,m}_ {P^{*}} \setminus \mathcal{V}^{X,m}_{P}$ are tested. Note, that if a nested approximation basis is used (i.e. $\Phibold_X(P) \subset \Phibold_X(P^{*})$), $\rbm(\ubm,\phibold)$ is a sub-vector of $\Rbm(\ubm,\phibold)$, so that only the latter needs to be computed. This fact is exploited in BoSSS where nested modal polynomial basis functions are employed.
	\paragraph{Constrained optimization problem}
	The objective function for the optimization problem $ f: \mathbb{R}^{N_U} \times \mathbb{R}^{N_s} \rightarrow \mathbb{R} $ is defined as
	\begin{eqnarray}
		f(\ubm,\varphibold) :=\frac{1}{2} \Rbm (\ubm,\varphibold)^T \Rbm (\ubm,\varphibold) = \frac{1}{2} \Vert \Rbm (\ubm,\varphibold) \Vert_2^2.
	\end{eqnarray}
	On its basis, the non-linear constrained optimization problem is formulated as: Find $(\ubm,\varphibold)$ such that
	\begin{eqnarray}
		\label{equation:Minimizing Problem}
		f(\ubm,\varphibold) \rightarrow \text{min!} \\ \text{subject to } \rbm(\ubm,\varphibold) =0. \nonumber
	\end{eqnarray}	
It is assumed that solutions $(\ubm,\varphibold)$ to this problem feature a discontinuity-aligned interface $\mathfrak{I}_\mathfrak{s}$, simultaneously satisfying the discretized system of conservation laws $\rbm (\ubm,\varphibold)=0$. 
	\section{Solving the constrained optimization problem} 
	In the following, details on a Sequential Quadratic Programming (SQP) method for solving the optimization problem \eqref{equation:Minimizing Problem} are provided. 
	\paragraph{Linear optimality system} The non-linear constrained optimization problem is addressed by introducing the corresponding Lagrange functional $\mathcal{L}: \mathbb{R}^{N_U} \times \mathbb{R}^{N_s} \times \mathbb{R}^{N_U} \rightarrow \mathbb{R} $ by \begin{eqnarray} 
		\mathcal{L}(\ubm, \varphibold, \lambdabold)=f(\ubm,\varphibold)-\lambdabold^{\top}\rbm(\ubm, \varphibold).
	\end{eqnarray} 
To solve the optimization problem, one typically aims to find a solution satisfying the following first order optimality condition: A pair $(\ubm^{*}, \varphibold^{*}) \in  \mathbb{R}^{N_U} \times \mathbb{R}^{N_s}$ is considered a first order solution if there exists a \( \lambdabold^* \in \mathbb{R}^{N_U} \) such that
	\begin{eqnarray}
		\label{equation:GradLagrange}
		\nabla \mathcal{L}\left(\ubm^*,\varphibold^*, \lambdabold^{*}\right)=0,
	\end{eqnarray}
	which in more detail writes as
	\begin{eqnarray}
		\begin{aligned}
			\nabla_u f\left( \ubm^*,\varphibold^*\right)-\frac{\partial \rbm}{\partial \ubm}\left( \ubm^*,\varphibold^*\right)^T \lambdabold^* &=0, \\
			\nabla_\varphi f\left( \ubm^*,\varphibold^*\right)-\frac{\partial \rbm}{\partial \varphibold}\left( \ubm^*,\varphibold^*\right)^T \lambdabold^* &=0, \\
			\textbf{\rbm}\left( \ubm^*,\varphibold^*\right) &=0.
		\end{aligned}
	\end{eqnarray}

	\paragraph{Newton's method} To solve \eqref{equation:GradLagrange} approximately, Newton's method is applied. Introducing the notation \begin{eqnarray}
		\zbm:= \begin{pmatrix}
			\ubm \\
			\varphibold
		\end{pmatrix} \in \mathbb{R}^{N_z}, \quad N_z=N_u + N_s,
	\end{eqnarray} Newton's method gives an iterative sequence $(\zbm_0,\lambdabold_0),(\zbm_1,\lambdabold_1), ... ,(\zbm_k,\lambdabold_k)$, converging to a solution $(\zbm^{*},\lambdabold^{*})=(\ubm^{*}, \varphibold^{*},\lambdabold^{*})$. In every Newton iteration, one typically solves the linear system 
\begin{equation}
	\Hbm_\mathcal{L}(\zbm_k,\lambdabold_k)\begin{pmatrix}
		\Delta\zbm_k \\
		\Delta \lambdabold_k
	\end{pmatrix} = - \nabla \mathcal{L}(\zbm_k,\lambdabold_k),
\end{equation}
for the Newton step $(\Delta \zbm_k~~\Delta \lambdabold_k)^T$ and it features the Hessian of $\mathcal{L}$ evaluated at the iterate $(\zbm_k,\lambdabold_k)$.
Here, the full Hessian $\Hbm_\mathcal{L}(\zbm,\lambdabold)$ of the Lagrangian writes as
	\begin{eqnarray}
		\Hbm_\mathcal{L}(\zbm,\lambdabold)=\begin{pmatrix}
			\Abm(\zbm,\lambdabold) & \Jbm_\rbm^T(\zbm) \\
			\Jbm_\rbm(\zbm) & 0
		\end{pmatrix}, \text{ where }\Abm(\zbm,\lambdabold) := \Hbm_f(\zbm) - \sum_{k=1}^{N_Z} \lambdabold_k \Hbm_{r_k}(\zbm).
	\end{eqnarray}
 Here, $\Hbm_f$ denotes the Hessian of the objective function $f$ (see \eqref{equation:Minimizing Problem}), $\Hbm_{r_k}$ denote Hessians for each residual components $r_k$ (see \eqref{eqn:localresidual}) and $\Jbm_\rbm$ denotes the Jacobian matrix of the full residual $\rbm$ (see \eqref{eqn:algebraic_residual}). 

\paragraph{Hessian approximation} Computing the Hessian $\Hbm_\mathcal{L}(\zbm,\lambdabold)$ is a very challenging task involving second-order derivatives that are not available in most CFD solvers. To circumvent this issue, we follow \cite{zahrImplicitShockTracking2020} and use the Levenberg-Marquardt method \cite{dennisNumericalMethodsUnconstrained1996}, replacing the Hessian with a partially regularized approximation
	\begin{eqnarray}
		\Abm(\zbm,\lambdabold)\approx \frac{\partial \Rbm}{\partial \zbm}(\zbm)^T \frac{\partial \Rbm}{\partial \zbm}(\zbm) + \gamma \begin{pmatrix}
			0 &0 \\
			0 & \Ibm
		\end{pmatrix}
	\end{eqnarray}
	where $\gamma \in \mathbb{R}^+$ is an regularization parameter and $\Ibm \in \Rbb^{N_s \times N_s}$ the identity matrix. We only regularize the components corresponding the DOFs of the level set (as it is assumed that the part belonging to the flow solution is invertible), so that the approximation of the Hessian becomes
	\begin{eqnarray}
		\Abm(\zbm,\lambdabold)\approx \Bbm(\zbm,\gamma):= \begin{pmatrix}
			\Bbm_{\ubm \ubm}(\zbm) & \Bbm_{\ubm \varphibold}(\zbm)  \\
			\Bbm_{\ubm \varphibold}(\zbm)^T & \Bbm_{\varphibold\varphibold}(\zbm,\gamma)
		\end{pmatrix},
	\end{eqnarray}
	where 
	\begin{align}
		\Bbm_{\ubm \ubm}(\zbm) &= \frac{\partial \Rbm}{\partial \ubm}(\zbm)^T \frac{\partial \Rbm}{\partial \ubm}(\zbm)
		\\
	\Bbm_{\ubm \varphibold}(\zbm)&=\frac{\partial \Rbm}{\partial \ubm}(\zbm)^T \frac{\partial \Rbm}{\partial \varphibold}(\zbm) \\
		\Bbm_{\varphibold \varphibold}(\zbm,\gamma)&= \frac{\partial \Rbm}{\partial \varphibold}(\zbm)^T \frac{\partial \Rbm}{\partial \varphibold} (\zbm)+ \gamma \Ibm.
	\end{align}
	Using this approach, one obtains the linear system
	\begin{eqnarray}
		\label{equation:linearsystem}
		\begin{pmatrix}
			\Bbm(\zbm_k)  & \Jbm^T_\rbm(\zbm_k) \\
			\Jbm_\rbm(\zbm_k) & 0
		\end{pmatrix} \begin{pmatrix}
			\Delta\zbm_k \\
			\Delta \lambdabold_k
		\end{pmatrix}
		= -\begin{pmatrix}
			\nabla f(\zbm_k)^T + \Jbm_\rbm^T(\zbm_k)\lambdabold_k \\
			\rbm(\zbm_k)
		\end{pmatrix}
	\end{eqnarray}
for the $k$-th iteration which needs to be solved to obtain the step $(\Delta \zbm_{k} ~~\Delta \lambdabold_k)^T$. 

\paragraph{Sequential Quadratic Programming method} Furthermore, the system can be reformulated to immediately solve for the Lagrange multiplier $\lambdabold_k$. The first line of the system reads as
	\begin{eqnarray}
		\Bbm(\zbm_k)\Delta\zbm_k +\Jbm_\rbm^T(\zbm_k)\underbrace{( \Delta \lambdabold_k +\lambdabold_k)}_{\lambdabold_{k+1}} =-\nabla f(\zbm_k)^T
	\end{eqnarray} and one can reformulate it to
	\begin{eqnarray}
		\label{equation:linearsystemDirect}
		\begin{pmatrix}
			\Bbm(\zbm_k)  & \Jbm^T_\rbm(\zbm_k) \\
			\Jbm_\rbm(\zbm_k) & 0
		\end{pmatrix} \begin{pmatrix}
			\Delta\zbm_k \\
			\lambdabold_{k+1}
		\end{pmatrix}
		= -\begin{pmatrix}
			\nabla f(\zbm_k)^T \\
			\rbm(\zbm_k)
		\end{pmatrix},
	\end{eqnarray}
	giving the typical update coming from a SQP method. Finally, after solving the linear system for $(\Delta \zbm_k ~~\lambdabold_k )^T$, the current iterate $\zbm_k$ is updated via
	\begin{eqnarray}
		\label{equation:solUpdate}
		\zbm_{k+1}  = 
		\zbm_{k} + \alpha_{k}\Delta\zbm_k,
	\end{eqnarray}
	where $\alpha_k \in [0,1) $ is a step length. To stabilize the method it has to be chosen carefully and it is computed via an inexact line search algorithm \cite{nocedalNumericalOptimization2006}. 
	
	\paragraph{Line search globalization}
	As commonly done in non-linear solvers, suitable step lengths $\alpha_k$ are computed via an inexact line search algorithm. This involves an iterative search for the maximal \begin{equation}
		\alpha_k \in \{ 1=\tau^0,\tau^{-1},\hdots , \alpha_{\text{min}}\}
	\end{equation} (for some $\tau \in (0,1)$), which satisfies the \textit{condition of sufficient decrease} for a merit function $\theta_k: \Rbb \rightarrow \Rbb$ and writes as
	\begin{eqnarray}
		\label{equation:sufficentDecreaseCondition}
		\theta_k(\alpha_k) \leq \theta_k(0) + \alpha_k \beta  \theta_k'(0),
	\end{eqnarray}
	where $\beta>0$ is a relaxation parameter. Following Nocedal et al. \cite{nocedalNumericalOptimization2006}, the non-smooth exact $l^1$-merit function	\begin{eqnarray}
		\theta_k(\alpha):= f(\zbm_k +\alpha \Delta\zbm_k) + \mu\Vert \rbm(\zbm_k + \alpha \Delta\zbm_k))\Vert_1
	\end{eqnarray}
is employed. Here,  $\mu(\zbm_k) := 2\Vert \lambdabold_{k+1} \Vert_\infty$ is controlling the relative weighting of the constraint and $\beta = 10^{-4}$, $\tau=0.5$ and $ \alpha_{\text{min}}=10^{-8}$ are chosen.
	
	\paragraph{Handling newborn cut-cells}
In the following discussion, we address how our method handles scenarios where the level set cuts a new cell. To ensure controlled interface movement and continuity with respect to our solution space, we enforce two key rules:
\begin{enumerate}
	\item Cells without a neighboring cut-cell shall not be cut in the current iteration. Any updates violating this condition are rejected, and the update step size $\alpha_k$ is adjusted accordingly. (i.e. \textit{make sure the interface moves only to neighbors})
	\item Solution extrapolation in newborn cut-cells: Before an uncut-cell is being cut, it belongs entirely to one of the sub-domains. However, when being cut by the interface, it splits into cut-cells that pertain to different sub-domains. In cases where cut-cells (referred to as "newborn" cut-cells) belong to sub-domains different from the original sub-domain, we implement a solution extrapolation strategy. (i.e. \textit{the solution values for newborn cut-cells are extrapolated from the largest neighboring cell within the same sub-domain, as determined by volume)}
\end{enumerate}
	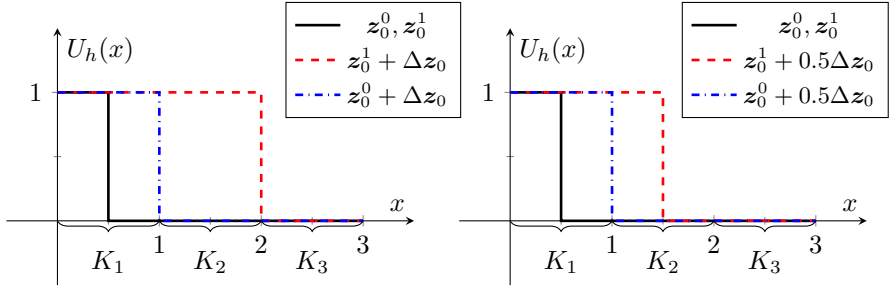
\begin{figure}[t]
	\begin{tikzpicture}
		\centering
		\begin{groupplot}[
			group style={group size=2 by 1, horizontal sep=0.6cm},
			xlabel={$x$},
			ylabel={$U_h(x)$},
			axis lines=middle,
			xmin=-0.5, xmax=3.5,
			ymin=-0.5, ymax=1.5,
			xtick={0,1,2,3},
			ytick={0,1},
			xticklabels={0,  1, 2, 3},
			yticklabels={0, 1},
			minor tick num=1,
			grid style={line width=0.1pt, draw=gray!30},
			major grid style={line width=0.2pt, draw=gray!60},
			width=7cm,
			height=5cm,
			]
			\nextgroupplot[
			legend style={at={(0.9,0.65)},anchor=south, font=\small},
			]
			\addplot[black, line width=1pt, mark size=3pt] coordinates {
				(0,1)
				(0.5,1)
				(0.5,0)
				(1,0)
				(2,0)
				(3,0)
			};
			\addplot[red, dashed,line width=1pt, mark size=3pt] coordinates {
				(0,1)
				(0.5,1)
				(1,1)
				(2,1)
				(2.0,0)
				(2.5,0)
				(3,0)
			};
			\addplot[blue, dashdotted,line width=1pt, mark size=3pt] coordinates {
				(0,1)
				(0.5,1)
				(1,1)
				(1,0)
				(2,0)
				(2.5,0)
				(3,0)
			};
			\addlegendentry{$\zbm_0^0,\zbm_0^1$}
			\addlegendentry{$\zbm_0^1 +\Delta \zbm_0$}
			\addlegendentry{$\zbm_0^0 + \Delta \zbm_0$}
			
			\draw[decorate,decoration={brace,amplitude=4pt,mirror}] (axis cs:0,0) -- (axis cs:1,0) node[midway,below=8pt]{$K_1$};
			\draw[decorate,decoration={brace,amplitude=4pt,mirror}] (axis cs:1,0) -- (axis cs:2,0) node[midway,below=8pt]{$K_2$};
			\draw[decorate,decoration={brace,amplitude=4pt,mirror}] (axis cs:2,0) -- (axis cs:3,0) node[midway,below=8pt]{$K_3$};
			\nextgroupplot[
			legend style={at={(0.8,0.65)},anchor=south, font=\small},
			]
			\addplot[black, line width=1pt, mark size=3pt] coordinates {
				(0,1)
				(0.5,1)
				(0.5,0)
				(1,0)
				(2,0)
				(3,0)
			};
			\addplot[red, dashed,line width=1pt, mark size=3pt] coordinates {
				(0,1)
				(0.5,1)
				(1,1)
				(1.5,1)
				(1.5,0)
				(2,0)
				(3,0)
			};
			\addplot[blue, dashdotted,line width=1pt, mark size=3pt] coordinates {
				(0,1)
				(0.5,1)
				(1,1)
				(1,0)
				(1.5,0)
				(2,0)
				(3,0)
			};
			\addlegendentry{$\zbm_0^0,\zbm_0^1$}
			\addlegendentry{$\zbm_0^1 + 0.5 \Delta \zbm_0$}
			\addlegendentry{$\zbm_0^0 + 0.5 \Delta \zbm_0$}
			
			\draw[decorate,decoration={brace,amplitude=4pt,mirror}] (axis cs:0,0) -- (axis cs:1,0) node[midway,below=8pt]{$K_1$};
			\draw[decorate,decoration={brace,amplitude=4pt,mirror}] (axis cs:1,0) -- (axis cs:2,0) node[midway,below=8pt]{$K_2$};
			\draw[decorate,decoration={brace,amplitude=4pt,mirror}] (axis cs:2,0) -- (axis cs:3,0) node[midway,below=8pt]{$K_3$};
		\end{groupplot}
	\end{tikzpicture}
	\caption{Exemplary application of an update $\Delta \zbm_0$ to an 1D XDG field $\zbm^l_0$ on a three cell grid $\{K_1,K_2,K_3\}$, illustrating how newborn cut-cells are handled. \textit{Left}: When applying $\zbm_0^l +  \Delta \zbm_0$, the interface moves from $x=0.5$ to $x=2.5$, which violates the first rule, i.e. the step length $\alpha=1$ is rejected. \textit{Right}: For $\zbm_0^l +  0.5\Delta \zbm_0$ the interface moves from $x=0.5$ to $x=1.5$, without violating the first rule. Solution is extrapolated in the new born cut-cell $K_2=[1,2]$ from the neighbor, i.e. choosing $\zbm_1=\zbm_0^1 + 0.5 \Delta \zbm_0$.}
	\label{fig:ls:movement}
\end{figure}
\paragraph{Example} To illustrate the application of the two rules concerning the handling of newborn cells, let us consider a simplified 1D example: Let $\Omega=[0,3]$ be a domain consisting of three cells $K_i=[i-1,i], i \in \nbrc{3}$ and divided into two sub-domains $\{\mathfrak{A},\mathfrak{B}\}$, where $\mathfrak{A} = \mathfrak{A}(\varphi):=  \left\{x\in \Omega \suchthat \varphi (x) <0\right\} ,\mathfrak{B} = \mathfrak{B}(\varphi):= \left\{ x\in \Omega \suchthat \varphi (x) >0\right\} $. Secondly, a globally parameterized level set $\varphi_0(x)=x+a$ having only a single DOF $a \in \mathbb{R}$, initially set to $a=-0.5$, is chosen. We then set up an element $U_h$ from the discrete $P=0$ XDG space by defining two coordinate vectors $\ubm^0,\ubm^1$ and a current iterate $	\zbm^l_0 $ by
	\begin{equation}
		\label{eq:ls_mov_coordvec}
		\ubm^l=(u^l_i)_{i \in \nbrc{6}}=(\underbrace{\underbrace{1}_{K_{1,\mathfrak{A}}},\underbrace{0}_{K_{1,\mathfrak{B}}}}_{\text{cut-cell}},\underbrace{\underbrace{l}_{K_{2,\mathfrak{A}}},\underbrace{0}_{K_{2,\mathfrak{B}}}}_{\text{uncut-cell}},\underbrace{\underbrace{0}_{K_{3,\mathfrak{A}}},\underbrace{0}_{K_{3,\mathfrak{B}}}}_{\text{uncut-cell}})^T, \quad	\zbm^l_0 =\begin{pmatrix}
			\ubm^1 \\
			a
			\end{pmatrix}
	\end{equation}
	for $l\in \{0,1\}$. We indicate the relationship between coefficients and cut-cells in \eqref{eq:ls_mov_coordvec} (assuming $P=0$ leads to one basis function per cut-cell, i.e. the indicator function supported by the respective cell). Note, that $K_{2,\mathfrak{A}} =K_{3,\mathfrak{A}}=\emptyset$ and thus $u^l_3,u^l_5$ do not affect the actual function value resulting in $U_h(x) =\ubm^0(x)\cdot \Phibold(x)=\ubm^1(x)\cdot \Phibold(x)$. Given an imaginary search direction $\Delta \zbm_0=(0,\hdots,0,-2)^T$ (only updating the level set DOF), we conduct a line search applying the rules from above: When adding the search direction $\Delta \zbm_0$ (depicted in Figure \ref{fig:ls:movement}) with $\alpha_0 =1$, one obtains
		\begin{equation}
	\zbm^l_1 = \zbm^l_0+\Delta \zbm_0=\begin{pmatrix}
		\ubm^l \\
		-2.5
	\end{pmatrix}.
		\end{equation}
	This implies that the updated level set becomes $\varphi_1(x)=x-2.5$, causing $K_3$ to become a cut-cell. In the resulting situation, $K_1$ is no longer a cut-cell and $K_2$ changes sub-domains.
	 However, this scenario violates the first rule, as the only neighbor ($K_2$) of $K_3$ was no cut-cell before the update. Consequently, the algorithm rejects $\alpha_0=1$ and proceeds to check $\tau\alpha_0=0.5$. 
	 Subsequently for $\alpha=0.5$, we obtain $\varphi_1(x)=x-1.5$, hence, $K_2$ becomes a cut-cell. Now, the first rule is not violated since $K_2$ had a neighboring cut-cell ($K_1$) before the update.
	 Further, the application of the second rule is illustrated in the following: Without solution extrapolation, one would obtain a different discrete state $\ubm^1\cdot \Phibold$, depending on the value of $l$, as $K_{2,\mathfrak{A}}$ is no longer a ghost cell. This resulting discrete state then writes
	\begin{equation}
		\ubm^1(x)\cdot \Phibold(x)=\begin{cases}
			1 &\text{for } x \leq 1 \\
			l &\text{for } 1< x \leq 1.5 \\
			0 &\text{else}.
		\end{cases}
	\end{equation}
To maintain continuity with respect to the solution space, we compute solution values in newborn cut-cells by extending the polynomial solution from the largest neighboring cell (in terms of volume) within the same sub-domain. For this example, $u^l_3=1$ is chosen when updating with $\alpha_0=0.5$, regardless of its previous value. 
	\paragraph{Adaptive choice of regularization parameter}
	The regularization parameter $\gamma \in \Rbb$ plays a crucial role in regulating the level-set-dependent part of our system's Hessian approximation. Its adaptability should be in accordance with the degree of change in the interface position and adaptivity is desirable. In simpler terms, if the level set update $\Delta \varphibold_k$ is perceived as excessively large ($\Vert \Delta\varphibold_k \Vert > \sigma_2 L$), more regularization is added, while if it is deemed excessively small ($\Vert \Delta\zbm_k \Vert < \sigma_1 L$), we aim for less regularization. This adaptability is achieved through the following approach:
	\begin{eqnarray}
		\gamma_{k+1} = \min \{\max \{\hat{\gamma}_{k+1}, \gamma_{min}\},\gamma_{max}\}, \quad \hat{\gamma}_{k+1}=
		\begin{cases}
			\kappa^{-1}\gamma_{k} &\text{if } \Vert \Delta \varphibold_k \Vert_2  < \sigma_1 L \\
			\kappa\gamma_{k} &\text{if } \Vert \Delta \varphibold_k \Vert_2 > \sigma_2 L\\
			\gamma_{k} &\text{else}.
		\end{cases}
	\end{eqnarray}
	Here, the values of $\sigma_1$, $\sigma_2$, and $L$ should be selected based on the length scales of the specific domain size. Additionally, the parameter $\kappa>1$ dictates the level of aggressiveness in the adaptation process.	For this study, we have opted for the following values: $L=1$, $\sigma_1=10^{-1}$, $\sigma_1=10^{-2}$, and $\kappa=1.5$.
	
	\paragraph{Termination/Change of basis criteria} 
	Determining when to terminate the method or to increase the solutions degree involves checking if the residual norms have dropped below a certain threshold, such as	$\Vert \rbm(\zbm_k) \Vert \leq \textrm{tol}$ and $\Vert \Rbm(\zbm_k) \Vert \leq \textrm{tol}$ Nevertheless, choosing a universal threshold value presents challenges. Setting it too low may result in the algorithm never terminating, while setting it too high could lead to premature termination, causing under-utilization of the high-order method. The challenge is addressed by continuing the method until the solution remains relatively stable over several specified steps $\termn$ (in this work $\termn=8$). To ensure robustness in this process, we introduce the residual-norm skyline, denoted as  $\text{sr}^c_{n} := \min_{j \leq n} | \cbm(\zbm_j) |$ ($\cbm\in \{ \rbm,\Rbm\}$), and the averaged reduction factor
	\begin{equation}
		\text{arf}^c_n := \frac{1}{\termn} \left(\sum_{k=n-\termn}^{n-1}\frac{ \text{sr}^c_{k} }{ \max({ \text{sr}_{k+1}, 10^{-100}}) }\right). 
	\end{equation} 
These metrics are defined for $n \geq \termn$ and applied to both residuals. The method is terminated if the following conditions are met for $\rbm$ and $\Rbm$:
	\begin{equation}
		\big( n \geq \termn \big) ~\wedge ~ \big(\textrm{sr}^c_n \leq 10^{-5} + 10^{-5} \Vert \ubm_n\Vert_2 \big) ~\wedge~ \big(\textrm{arf}^c_n < 1.001\big).
	\end{equation}
	This choice ensures that the non-linear system is solved with high accuracy within the constraints of floating-point precision. Additionally, the skyline approach enhances robustness against oscillations near the lower limit.
	\section{Level set parametrization}
	At this point, we aim for a reduction of DOFs for the interface representation. It is of interest for the assembly of the linear system \eqref{equation:linearsystem} as it requires the computation of derivatives $\pder{\rbm}{\varphibold}, 	\pder{\Rbm}{\varphibold}$ of the residuals $\rbm,\Rbm$ with respect to coefficients $\varphibold$ defining the level set. This amounts to differentiation of integrals on implicitly defined surfaces and we revert to central finite differences for their computation, i.e.
	\begin{equation}
		\pder{\rbm}{\varphi_l} \approx \frac{\rbm(\ubm,\hdots ,\varphi_{l-1},\varphi_l +\epsilon, \varphi_{l+1},\hdots) -\rbm(\ubm,\hdots ,\varphi_{l-1},\varphi_l -\epsilon, \varphi_{l+1},\hdots)}{2 \epsilon},
	\end{equation} 
where $\varphi_l$ denotes one exemplary component of $\varphibold$ and $\epsilon>0$ (typically $\epsilon=10^{-8}$). Consequently, for each perturbed level set, new quadrature rules for cut-cells must be computed. BoSSS accomplishes this task using an algorithm developed by Saye \cite{sayeHighOrderQuadratureMethods2015} based on non-linear root-finding. This process incurs significant computational costs, rendering the differentiation method impractical when the level set representation involves a high number of DOFs which essentially amounts to evaluating the residuals for numerous slightly different level sets and interface positions. Typically, in a XDG method, interfaces are represented implicitly by level sets residing in the DG space $\mathcal{V}_{P_s}$ sharing the same grid as the corresponding solution. This introduces many DOFs because of the high dimension of $\mathcal{V}_{P_s}$ and redundancy (e.g., $\varphi$ and $a \varphi$ define the same interface for any constant $a\in \Rbb$). Also, we noticed that allowing discontinuities/kinks in the interface across cell boundaries has negative effects on the success of the optimizer, leading to more stagnation. 

\paragraph{Spline level sets} To address these challenges, we simplify the level set representation by directly parametrizing the interface using $C^1$-continuous splines: Given a set of interpolation points $y_0 <y_1 <\hdots<y_{N_S}$ (corresponding to the $y$ coordinates of the vertices of the Cartesian cells of the grid), along with associated values  $\{x_0,x_1,...,x_{N_S}\}$ and derivatives $\{c_0,c_1,...,c_{N_S}\}$ at these points, one can define a cubic spline $S:[y_0,y_{N_S}] \rightarrow \mathbb{R}$ by
	\begin{equation}
		S(y_i)=x_i, ~S'(y_i)=c_i, ~S_{\vert_{[y_i,y_{i+1}]}} \in \mathbb{P}_3([y_i,y_{i+1}]).
	\end{equation} 
	The resulting spline is a $C^1$-function defined on the interval $[y_0, y_S]$. To transform it into a discrete level set, one can project the function
	\begin{equation}
		\label{equation:spline}
		\varphi_S(x,y) = x-S(y)
	\end{equation}
	onto the DG space $\mathcal{V}_P$. This representation has a fixed number of two DOFs per grid point in the $y$ direction, significantly fewer than a DG field of polynomial order three would have. Furthermore, it explicitly represents the shock front. While this specific spline representation can be extended to any polynomial order (for instance, a continuous linear spline can be defined similarly by omitting the derivative information), this work primarily utilizes cubic and linear splines.
	
	It is important to acknowledge that this representation is a significant simplification, limited to representing shocks that correspond to the graph of a 1D height function. However, for the purposes of a proof of concept and method development, this choice is considered valid, as there are numerous instances of flows with shocks that can be adequately represented by a spline. In future work, we aim to develop a more general interface representation together with an efficient algorithm to obtain the corresponding derivatives.
	\section{Robustness measures}
Numerical experiments have indicated that the method can exhibit convergence issues or get stuck at undesirable local minima in certain scenarios. These issues are more prominent when dealing with small cut-cells and when the polynomial degree $P$ is chosen greater than zero in cases where the current iterate is significantly distant from the exact solution. Consequently, we introduce additional stabilization measures to enhance the method's robustness, each addressing distinct challenges that may arise.
	\subsection{Agglomeration}
	One of the identified causes of stagnation is the presence of small cut-cells. Since the background mesh remains Cartesian and does not adapt to the shock front's position, the optimal interface location may lead to the creation of small cut-cells. Additionally, during the optimization process, intermediate iterations may also exhibit small cut-cells, due to arbitrary interface positions. These small cut-cells pose an issue, introducing an imbalance in the solution vector, as their coefficients inversely scale with the cell size. This in turn magnifies the system's condition number, leading to a destabilization of the method \cite{deprenterStabilityConditioningImmersed2022}.
	
	To address this issue, a technique known as \textit{cell agglomeration} is employed. As the name implies, this approach involves agglomerating cut-cells $K_{j,\mathfrak{s}}$ for which the ratio
	\begin{equation}
		\label{eqn:aggloCondition}
		\frac{\vert K_{j,\mathfrak{s}}\vert}{\vert K_j \vert} \leq \text{agg}_{\text{thrsh}}
	\end{equation} falls below a predefined threshold $\text{agg}_{\text{thrsh}} \in (0.1)$. These cut-cells are merged with their largest neighbor within their respective sub-domain, while elements on the opposite side of a discontinuity remain unaffected. Consequently, this process transforms the cut-cell grid $\mathcal{K}^X_h$ into an agglomerated cut-cell grid $\text{Agg}(\mathcal{K}^X_h,\Acal)$, represented by an agglomeration map $\Acal$, a subset of all edges shared by two cut-cells belonging to the same domain, i.e.
	\begin{equation}
		 \Acal \subset  \left\{(K_{i,s}, K_{j,s}) \suchthat \overline{K}_{i,s} \cap \overline{K}_{j,s}  \neq \emptyset, (i,j) \in [J]^2, s \in \text{SD} \right\}.
	\end{equation}
 Here, the relation $(K_{i,s}, K_{j,s})\in \Acal$ indicates that $K_{i,s}$ is agglomerated into $K_{j,s}$ and $\Acal$ is determined by the condition \eqref{eqn:aggloCondition}. In the context of BoSSS, \textit{cell agglomeration} has already proven successful in various applications, with the development of a dedicated multigrid algorithm built around it and we refer the reader to \cite{kummerBoSSSPackageMultigrid2021} for details.
	
	In the specific context of this optimization method, \textit{cell agglomeration} is applied as follows:
	
	\begin{enumerate}
		\item For each iteration step, residuals  $\rbm(\zbm_k), \Rbm(\zbm_k)$ and sensitivities $\Jbm_\rbm(\zbm_k)$, $\Jbm_R(\zbm_k)$ are computed on the non-agglomerated mesh $\Kcal^X_h$.
		\item The agglomeration map $\Acal_k$ is computed based on \eqref{eqn:aggloCondition} and depending on $\varphibold_k$.
		\item The agglomerated state is computed $\zbm^\text{agg}_k$ from $\zbm_k$, necessitating a basis change from the non-agglomerated basis of $\Vcal^{X,m}(\mathcal{K}^X_h)$ to the agglomerated basis of $\Vcal^{X,m} (\text{Agg}(\mathcal{K}^X_h,\Acal_k))$.
		\item The agglomerated residuals $\rbm^\text{agg}(\zbm^\text{agg}_k), \Rbm^\text{agg}(\zbm^\text{agg}_k)$ and sensitivities $\Jbm^\text{agg}_\rbm(\zbm^\text{agg}_k)$, $\Jbm^\text{agg}_\Rbm(\zbm^\text{agg}_k)$ are computed from their non-agglomerated counterparts by the same basis change.
		\item The agglomerated version of the linear system \eqref{equation:linearsystemDirect} is assembled and solved to determine the agglomerated inexact Newton step $\Delta\zbm_k^{agg}$.
		\item During each line search iteration $ \zbm^\text{agg}_k + \alpha \Delta\zbm^\text{agg}_k$, agglomerated residuals $\rbm^\text{agg}( \zbm^{agg}_k + \alpha \Delta\zbm^\text{agg}_k), \Rbm^\text{agg}( \zbm^\text{agg}_k + \alpha \Delta\zbm^\text{agg}_k)$ are computed when evaluating the merit function.
		\item After finding a suitable step length $\alpha_k$, the agglomerated iterate $\zbm^\text{agg}_{k+1}=\zbm^\text{agg}_k + \alpha_k \Delta\zbm^\text{agg}_k$ is finally extrapolated to the non-agglomerated mesh to obtain the non-agglomerated iterate $\zbm_{k+1}$.
	\end{enumerate}
	So far, the use of agglomeration within the XDG shock tracking method has shown a positive impact on the quality of search directions $\Delta \zbm_k$, primarily acting as a preconditioner for the linear system. 
	
	\subsection{$P$-continuation strategy}
	In our numerical experiments we observed that when significant level set movement is required, particularly in cases where the initial guess for the level set is inaccurate, a beneficial strategy involves initially restricting the polynomial degree to $P=0$ and gradually increasing it. This gradual increment occurs when the solver fails to compute substantial updates to the solutions, and the minimum required iterations $\{\text{min}_0,\text{min}_1,\hdots \}$ for the respective degrees are exceeded. In our work, we selected $\{\text{min}_0=30,\text{min}_1=30,\text{min}_2=10,\text{min}_3=10 \}$ for these minimum iterations. The decision to increase the polynomial degree is determined by checking the termination criteria described in Section 5. 
	
		\subsection{Solution re-initialization in oscillatory cells}
		In accordance with Huang \& Zahr \cite{huangRobustHighorderImplicit2022a}, we introduce a re-initialization procedure to enhance robustness. Similar to the authors' findings, we observed the emergence of non-physical oscillations in the XDG solution during intermediate steps when employing polynomial degrees greater than zero. These oscillations subsequently lead to poor search directions in subsequent iterations and result in very small step sizes $\alpha_k$.
		
		To address this issue, we employ the following re-initialization procedure: It identifies oscillatory cells using a variant of the well-established Persson-Peraire shock sensor. In these cells, the XDG solution is reset to a $P=0$ field by averaging over a patch of elements in the cell's vicinity.
		
		\paragraph{Determining oscillatory cells} For a cell $K=K_{j,\mathfrak{s}} \in \mathcal{K}^X_h(\varphi) =\mathcal{K}^X_h$ (omitting the level set dependency for readability) and an XDG solution $U=U_h \in \Vcal^{X,m}_P$ of polynomial degree $P>0$, the Persson-Peraire shock sensor $ S: \mathcal{K}^X_h \times \mathcal{V}^{X,m}_P \rightarrow \mathbb{R}$ \cite{perssonSubCellShockCapturing2006} for the first component $U_1\in \Vcal^{X,1}_P$ is defined as 
		\begin{equation}
			\label{eq:pers_sensor}
				S(K,U) := \log\left(\frac{\Vert U_1 - \pi_{P-1}(U_1) \Vert_{L^2(K)}}{\Vert U_1 \Vert_{L^2(K)}}\right),
			\end{equation}
		where $\pi_{P-1} : \mathcal{V}^{X,1}_P \rightarrow \mathcal{V}^{X,1}_{P-1}$ denotes the orthogonal projection onto the XDG space $\mathcal{V}^{X,1}_{P-1}$ with polynomial degree $P-1$.
		Using this, one can define the set of all oscillatory cells by
		\begin{equation}
				\mathcal{K}_\text{Osc}(U) :=\left\{ K \in \mathcal{K}^X_h \suchthat S(K ,U) > \epsilon_1\right\} ,
			\end{equation}
		where $\epsilon_1>0$ (typically we choose $\epsilon_1=-0.2$) is a threshold chosen depending on the problem. Huang \& Zahr \cite{huangRobustHighorderImplicit2022a} propose not only re-initializing the oscillatory cells but also their neighbors. This results in the final set of cells to be re-initialized
		\begin{equation}
				\mathcal{K}_\text{ReInit}(U) := \mathcal{K}_\text{Osc}(U) \cup \mathcal{N}(\mathcal{K}_\text{Osc}(U)),
			\end{equation}
		where $\mathcal{N}(\mathcal{K}_\text{Osc}(U))$ represents the subset of cells sharing a face with any cell in $\mathcal{K}_\text{Osc}(U)$. In cases of excessive line search iterations, typically occurring beyond a predefined threshold (suggested by the authors as 5 iterations), the set is redefined as:
		\[
		\mathcal{K}_\text{ReInit}(U) :=  \left\{K \in \mathcal{K}^X_h  \suchthat S(K ,U) > \epsilon_2 \max_{K \in \mathcal{K}^X_h }S(K ,U)\right\},
		\]
		where we adopt the authors' choice of $\epsilon_2 = 10^{-2}$.
		\paragraph{Shock-aware solution re-initialization} The XDG solution is reset on the cells $K \in \mathcal{K}_\text{ReInit}(U) $ to a constant value, which is the average of the current solution $U$ over a subset of its neighboring cells $\tilde{\mathcal{N}}(K,U)$ defined in the following. Instead of considering all neighboring cells, Huang \& Zahr choose to only include those on the same side of the discontinuity/shock. As to determine these cells, the average jump function $ a_J: \mathcal{K}^X_h \times \mathcal{K}^X_h\times \mathcal{V}^{X,m}_P \rightarrow \mathbb{R}$ is defined as
		\begin{equation}
				a_J(K,K',U) \mapsto \frac{1}{\vert \partial K \cap \partial K'\vert} \int_{\vert \partial K \cap \partial K'\vert} \jump{U_1} dS,
			\end{equation}
		where $\jump{U_1}$ represents the jump in the $U_1$ component across the common face of cells $K$ and $K'$. Note that, $a_J$ is also helpful in selecting cells on only one side of the shock if it does not coincide with the interface $\mathfrak{I}_s$. Consequently, the modified neighbor set
		\begin{equation}
				\tilde{\mathcal{N}}(K,U) := \left\{K' \in \mathcal{N}(K) \suchthat \vert a(K,K',U) \vert\leq \epsilon_3\right\} ,
			\end{equation}
		is defined with $\epsilon_3=10^{-2}$. On this patch, the solution is reset to constant values, which can be expressed by
		\begin{equation}
				U^{\text{ReInit}}_i\vert_{K}:= \begin{cases}
						\frac{1}{\vert K \cup \tilde{\mathcal{N}}(K,U) \vert} \int_{\vert K \cup \tilde{\mathcal{N}}(K,U) \vert} U_i dV 
						&\text{if } K \in \mathcal{K}_\text{ReInit}(U)
						\\
						U_i &\text{else}
					\end{cases} \text{ for } i=1,\hdots,m.
			\end{equation}
		Lastly, re-initialization is prohibited if the current iterate is already close to satisfying the constraints, i.e
		\begin{equation}
			\label{eq:reInit_eps4}
				\Vert \rbm(\zbm_k) \Vert \leq \epsilon_4,
			\end{equation}
		where $\epsilon_4=10^{-2}$. Also, after a certain number of iterations $M=30$, re-initialization is no longer applied to ensure asymptotic convergence of the method. In the setting of $P$-continuation, this number is increased after each change of the polynomial degree of the solution to allow for re-initialization for the higher degrees.
		
		The complete algorithm for the solution re-initialization on oscillatory cells is shown in Algorithm \ref{alg:reInit}. Depending whether it is used due to excessive line search ($\textit{isEL}=\textit{true}$) or not ($\textit{isEL}=\textit{false}$) the solution may be re-initialized on a different set of cells $\mathcal{K}_\text{ReInit}(U)$.
		\begin{algorithm}
			\caption{Solution re-initialization on oscillatory cells}
			\label{alg:reInit}
			\begin{algorithmic}[1]
				\REQUIRE XDG solution $U$, grid $\mathcal{K}_h^X$, parameters $\epsilon_1,\epsilon_2,\epsilon_3$ and logical $\textit{isEL}$ (indicating whether solution re-initialization is applied due to excessive line search)
				\ENSURE XDG solution $U^{\text{ReInit}}$ re-initialized on oscillatory cells 
				\STATE \textbf{Compute max shock sensor:} $S_{\text{max}}:=\max_{K \in \mathcal{K}^X_h } S(K,U)$
				\STATE \textbf{Compute set of cells to be re-initialized $\mathcal{K}_\text{Osc}(U)$:}
				\STATE $\mathcal{K}_\text{Osc}(U) := \emptyset $ (initialize)
				\FOR{$K \in \mathcal{K}^X_h$}
					\STATE Compute shock sensor $S(K,U)$
					\IF{$\textit{isEL}$ and $S(K,U)>\epsilon_2 S_{\text{max}}$} 
						\STATE $\mathcal{K}_\text{Osc}(U):=\mathcal{K}_\text{Osc}(U) \cup K$ (add cell)
					\ELSIF{$S(K,U)>\epsilon_1$} 
						\STATE $\mathcal{K}_\text{Osc}(U):=\mathcal{K}_\text{Osc}(U) \cup K$ (add cell)
					\ENDIF
				\ENDFOR
				\STATE  $\mathcal{K}_\text{ReInit}(U) := \mathcal{K}_\text{Osc}(U) \cup \mathcal{N}(\mathcal{K}_\text{Osc}(U))$ (add neighbors)
				\STATE \textbf{Re-initialize solution:}
				\STATE $U^{\text{ReInit}}:=U$ (initialize)
				\FOR{$K \in \mathcal{K}_\text{ReInit}(U)$}
					\STATE $\tilde{\mathcal{N}}(K,U):=\emptyset$ (initialize)
					\FOR{$K' \in \mathcal{N}(K)$}
						\IF{$a_J(K,K',U) \leq \epsilon_3$}
							\STATE $\mathcal{N}(K,U):=\mathcal{N}(K,U) \cup K'$
						\ENDIF
					\ENDFOR
					\STATE Set $U^{\text{ReInit}}\vert_{K}:=	\frac{1}{\vert K \cup \tilde{\mathcal{N}}(K,U) \vert} \int_{\vert K \cup \tilde{\mathcal{N}}(K,U) \vert} U dV$ 
				\ENDFOR
				\STATE Set $U=U^{\text{ReInit}}$
			\end{algorithmic}
		\end{algorithm}
	\section{Results}
	In this section, we present results for three distinct conservation laws characterized by a solution discontinuity, the  1D space-time advection equation, the 1D space-time Burgers equation and the steady 2D inviscid Euler equations. The cases are presented in ascending order of complexity.
	\subsection{Linear advection equation}
	\begin{figure}[h]
		\raisebox{0.5cm}{\begin{subfigure}{0.08\textwidth}
				\centering
				\begin{tikzpicture}
					\begin{axis}[
						title={$c$},
						xlabel={},
						ylabel={},
						xmin=0, xmax=1,
						ymin=-3, ymax=4,
						xtick=\empty,
						ytick={-3,0.5,4},
						yticklabels={0,0.5,1},
						grid=both,
						width=4.8cm,  
						axis equal image
						]
						\addplot graphics [xmin=0,xmax=1,ymin=-3,ymax=4] {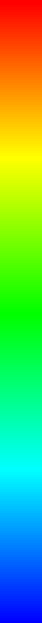};
					\end{axis}
				\end{tikzpicture}
		\end{subfigure}}
		\quad
		\begin{subfigure}{.33\textwidth}
			\centering
			\begin{tikzpicture}
				\begin{axis}[
					title=Iteration 0,
					ylabel={$t$},
					xlabel={$x$},
					xlabel style={yshift=5pt}, 
					ylabel style={yshift=-10pt}, 
					xmin=0, xmax=768,
					ymin=0, ymax=819,
					xtick={0,384,768},
					xticklabels={0,0.5,1},
					ytick={0,409,819},
					yticklabels={0,0.5,1},
					axis equal image,
					grid=both,
					width=1.5\linewidth, 
					]
					\addplot graphics [xmin=0,xmax=768,ymin=0,ymax=819] {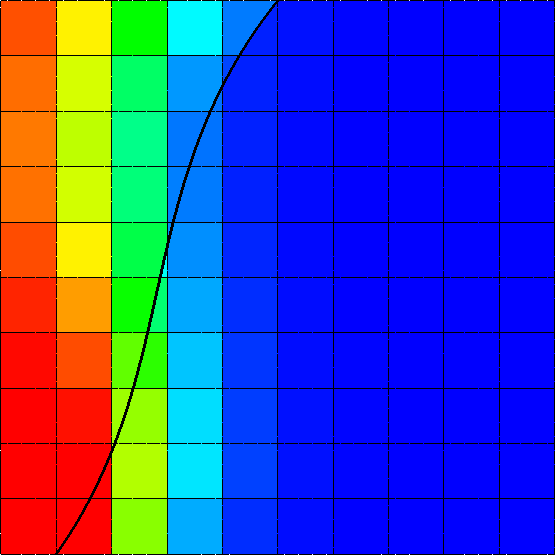};
				\end{axis}
			\end{tikzpicture}
		\end{subfigure}%
		\quad \quad
		\begin{subfigure}{.33\textwidth}
			\centering
			\begin{tikzpicture}
				\begin{axis}[
					title=Iteration 5,
					xlabel={$x$},
					xlabel style={yshift=5pt},
					xmin=0, xmax=768,
					ymin=0, ymax=819,
					xtick={0,384,768},
					xticklabels={0,0.5,1},
					ytick=\empty,
					axis equal image,
					grid=both,
					width=1.5\linewidth, 
					]
					\addplot graphics [xmin=0,xmax=768,ymin=0,ymax=819] {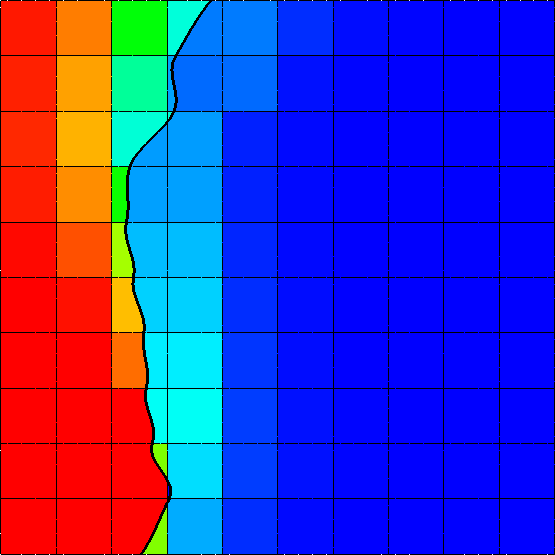};
				\end{axis}
			\end{tikzpicture}
		\end{subfigure}
		\\
		\begin{subfigure}{.33\textwidth}
			\centering
			\begin{tikzpicture}
				\begin{axis}[
					title=Iteration 10,
					ylabel={$t$},
					ylabel style={yshift=-10pt},
					xlabel={$x$},
					xlabel style={yshift=5pt},
					xmin=0, xmax=768,
					ymin=0, ymax=819,
					xtick={0,384,768},
					xticklabels={0,0.5,1},
					ytick={0,409,819},
					yticklabels={0,0.5,1},
					axis equal image,
					grid=both,
					width=1.5\linewidth, 
					]
					\addplot graphics [xmin=0,xmax=768,ymin=0,ymax=819] {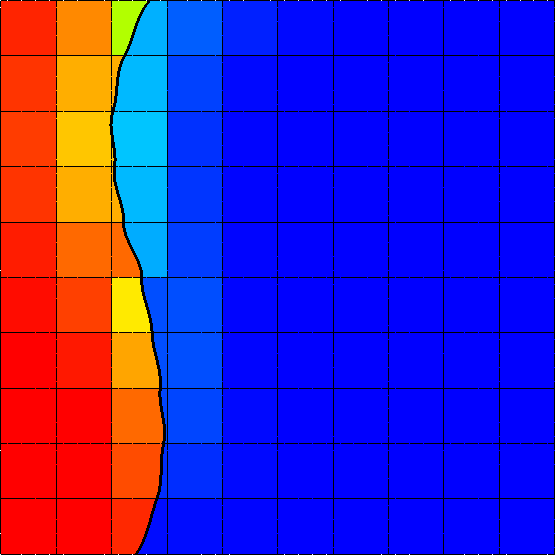};
				\end{axis}
			\end{tikzpicture}
		\end{subfigure}%
		\quad \quad
		\begin{subfigure}{.33\textwidth}
			\centering
			\begin{tikzpicture}
				\begin{axis}[
					title=Iteration 15,
					xlabel={$x$},
					xlabel style={yshift=5pt},
					xmin=0, xmax=768,
					ymin=0, ymax=819,
					xtick={0,384,768},
					xticklabels={0,0.5,1},
					ytick=\empty,
					axis equal image,
					grid=both,
					width=1.5\linewidth, 
					]
					\addplot graphics [xmin=0,xmax=768,ymin=0,ymax=819] {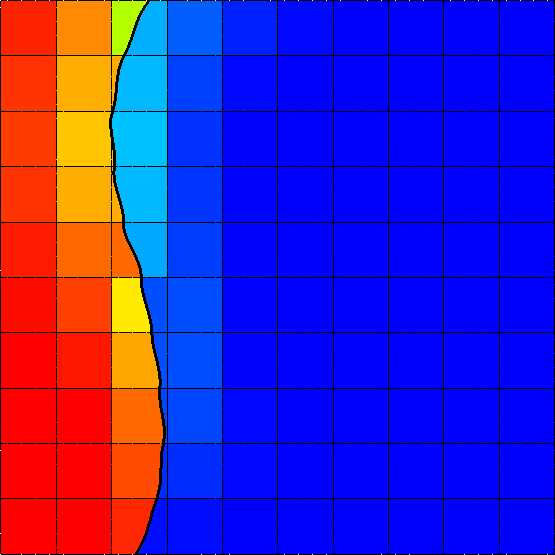};
				\end{axis}
			\end{tikzpicture}
		\end{subfigure}%
		\begin{subfigure}{.33\textwidth}
			\centering
			\begin{tikzpicture}
				\begin{axis}[
					title=Iteration 50,
					xlabel={$x$},
					xlabel style={yshift=5pt},
					xmin=0, xmax=768,
					ymin=0, ymax=819,
					xtick={0,384,768},
					xticklabels={0,0.5,1},
					ytick=\empty,
					axis equal image,
					grid=both,
					width=1.5\linewidth, 
					]
					\addplot graphics [xmin=0,xmax=768,ymin=0,ymax=819] {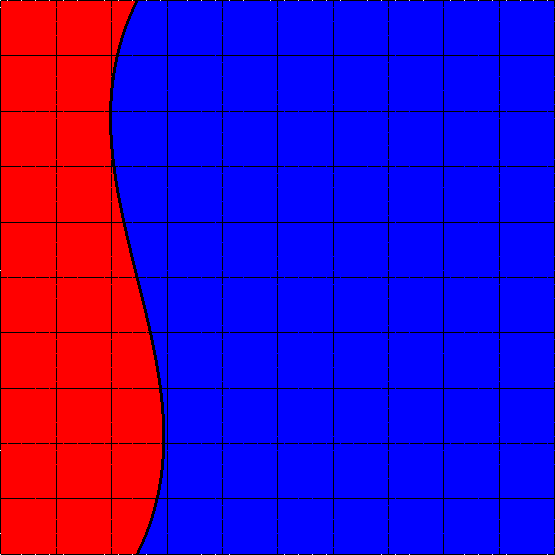};
				\end{axis}
			\end{tikzpicture}
		\end{subfigure}
		
		\caption{Plots of concentration $c$ for selected XDG shock tracking iterations for the space-time advection test case ($P=0$) for $k=0,5,10,15,50$. The curved polynomial discontinuity is tracked by the SQP solver starting from an initial guess for the shock interface $\mathfrak{I}_s$ (\textit{thick black line}), implicitly defined by a cubic spline level set ($P_s=3$) and which is far away from the correct position. Simultaneously, the solution of the scalar conservation law is obtained and the cubic spline level set is successfully aligned with the discontinuity.}
		\label{fig:SAIDT}
	\end{figure}
	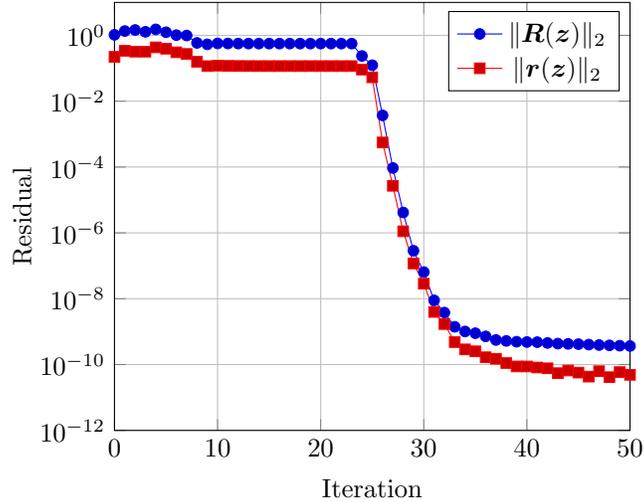
\begin{figure}[h]
		\centering
		\begin{tikzpicture}
			\begin{semilogyaxis}[            xlabel={Iteration},            ylabel={Residual},            xmin=0, xmax=50,            ymin=1e-12, ymax=10,            xtick={0,10,20,30,40,50},            ytick={1,1e-2,1e-4,1e-6,1e-8,1e-10,1e-12},            xticklabel style={/pgf/number format/.cd,fixed,precision=0},            yticklabel style={/pgf/number format/.cd,fixed,precision=0,sci},            grid=both,            ]
				\addplot table[x index=0,y index=1] {CurvedShockAdvection/R1norms.txt};
				\addlegendentry{$\Vert \Rbm(\zbm) \Vert_2$}
				\addplot table[x index=0,y index=1] {CurvedShockAdvection/R0norms.txt};
				\addlegendentry{$\Vert \rbm(\zbm) \Vert_2$}
			\end{semilogyaxis}
		\end{tikzpicture}
		\caption{Optimization history of the XDG shock tracking for the space-time advection test case ($P=0$). Both residual norms $\Vert \Rbm \Vert_2, \Vert\rbm \Vert_2$ converge rapidly after 30 iterations as the SQP solver successfully tracks the polynomial discontinuity.}
		\label{fig:OptHisSAIDT}
	\end{figure} 
	First, the method is showcased for the space-time formulation of the linear advection equation which writes
	\begin{eqnarray}
		\label{eq:st_adv}
		\frac{\partial c}{\partial t} + a(t)\frac{\partial c}{\partial x} =0 \quad \text{in }\Omega=[0,1]\times[0,1],
	\end{eqnarray}
	where
	\( a:[0,1] \rightarrow \mathbb{R} \) represents the time-dependent velocity and $c:\Omega \rightarrow \mathbb{R}$ a concentration. Here, the corresponding physical flux is $\fbm(c)=(a(t)c~ c)$.
	This equation is augmented with a discontinuity in the initial value
	\begin{eqnarray} \label{eq:advec_initial} 
		c(x,0)=H \left( \frac{1}{4} -x\right), \end{eqnarray} 
	where $H$ denotes the Heaviside function. Consequently, this problem involves the advection of a discontinuity in time. We impose Dirichlet boundary conditions using the exact solution
	\begin{eqnarray} c_\text{Ex}(x,t) = H(s(t)-x).\end{eqnarray} 
	Here, the height function $s: \Rbb \rightarrow \Rbb$ describing the discontinuity is defined as
	\begin{eqnarray}
		s(t) = \frac{1}{4} +\int_0^{t} a(t') dt'.
	\end{eqnarray}
	To discretize the conservation law \eqref{eq:st_adv}, we employ a Cartesian $10 \times 10$ background mesh. As the numerical flux, we select the upwind flux
	\begin{eqnarray}
		\hat{F}(c^+,c^-,\nbm_\Gamma)= \begin{cases}c^- (a(t),1)^T \cdot \nbm_\Gamma \quad \text{if } (a(t),1)^T \cdot \nbm_\Gamma <0 \\
			c^+ (a(t),1)^T\cdot \nbm_\Gamma \quad \text{if } (a(t),1)^T\cdot \nbm_\Gamma \geq 0 \end{cases}.
	\end{eqnarray}

	 Further, we aim to construct a problem configuration which serves to assess the shock tracking method's capability to align a curved interface, even when the initial guess for its position lacks sub-cell accuracy. According to this premise, we choose $a(t)$ as 
	\begin{eqnarray}
		a(t)=3t^2-3t+\frac{1}{2} .
	\end{eqnarray}
	The exact solution of this problem is piecewise constant and can be exactly represented in an XDG space of order $P=0$ using the level set of order  $P_s=3$
	\begin{eqnarray}
		\varphi_s(x,t)=x-s(t)=x-t^3+\frac{3}{2}t^2-\frac{1}{2}t-\frac{1}{4}.
	\end{eqnarray} Therefore, we set the polynomial degree to $P=0$ and use a cubic spline level set $\varphi_s$. For the initial guess we set
	\begin{eqnarray}
		\varphi_s(x,t)=x-s(t)=x-\frac{7}{10}t^3+t^2-\frac{7}{10}t-\frac{1}{10},
	\end{eqnarray} and project the exact solution for the state, as depicted in Figure \ref{fig:SAIDT}. The initial guess for the level set significantly deviates from the actual solution in terms of position, while the curvature is close. 

Finally, we apply the XDG shock tracking method to the described configuration. After 30 iterations, the method converges, perfectly aligning the discontinuity and reducing the residual to $10^{-10}$. Figure \ref{fig:SAIDT} displays plots of selected iterations of the SQP solver, while Figure \ref{fig:OptHisSAIDT} illustrates the history of both residual norms.
	\subsection{Inviscid Burgers equation}
	\begin{figure}[t]
		\raisebox{0.5cm}{\begin{subfigure}{0.08\textwidth}
				\centering
				\begin{tikzpicture}
					\begin{axis}[
						title={$c$},
						xlabel={},
						ylabel={},
						xmin=0, xmax=1,
						ymin=-3, ymax=4,
						xtick=\empty,
						ytick={-3,0.5,4},
						yticklabels={0.25,0.5,0.75},
						grid=both,
						width=4.8cm,  
						axis equal image
						]
						\addplot graphics [xmin=0,xmax=1,ymin=-3,ymax=4] {legend_trimmed.png};
					\end{axis}
				\end{tikzpicture}
		\end{subfigure}}
		\quad
		\begin{subfigure}{.33\textwidth}
			\centering
			\begin{tikzpicture}
				\begin{axis}[
					title=Iteration 0,
					ylabel={$t$},
					xlabel={$x$},
					xlabel style={yshift=5pt}, 
					ylabel style={yshift=-10pt}, 
					xmin=0, xmax=768,
					ymin=0, ymax=819,
					xtick={0,384,768},
					xticklabels={0,0.5,1},
					ytick={0,409,819},
					yticklabels={0,0.5,1},
					axis equal image,
					grid=both,
					width=1.5\linewidth, 
					]
					\addplot graphics [xmin=0,xmax=768,ymin=0,ymax=819] {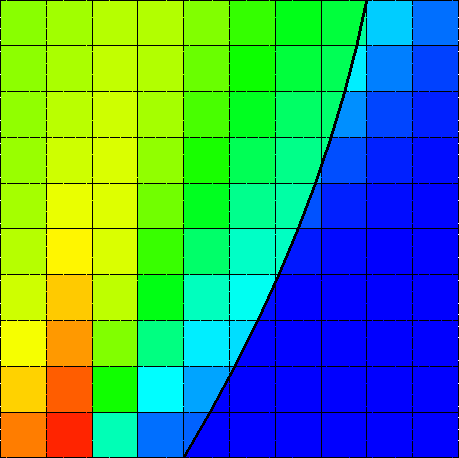};
				\end{axis}
			\end{tikzpicture}
		\end{subfigure}%
		\quad \quad
		\begin{subfigure}{.33\textwidth}
			\centering
			\begin{tikzpicture}
				\begin{axis}[
					title=Iteration 2,
					xlabel={$x$},
					xlabel style={yshift=5pt},
					xmin=0, xmax=768,
					ymin=0, ymax=819,
					xtick={0,384,768},
					xticklabels={0,0.5,1},
					ytick=\empty,
					axis equal image,
					grid=both,
					width=1.5\linewidth, 
					]
					\addplot graphics [xmin=0,xmax=768,ymin=0,ymax=819] {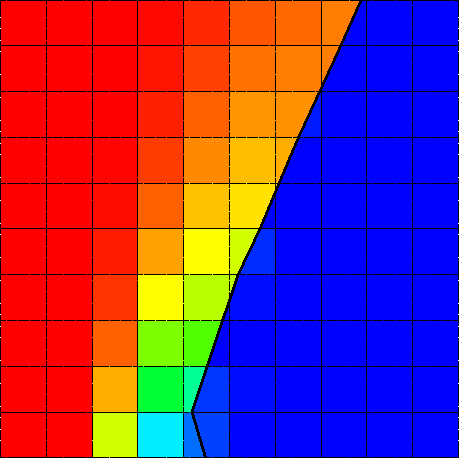};
				\end{axis}
			\end{tikzpicture}
		\end{subfigure}
		\\
		\begin{subfigure}{.33\textwidth}
			\centering
			\begin{tikzpicture}
				\begin{axis}[
					title=Iteration 5,
					ylabel={$t$},
					ylabel style={yshift=-10pt},
					xlabel={$x$},
					xlabel style={yshift=5pt},
					xmin=0, xmax=768,
					ymin=0, ymax=819,
					xtick={0,384,768},
					xticklabels={0,0.5,1},
					ytick={0,409,819},
					yticklabels={0,0.5,1},
					axis equal image,
					grid=both,
					width=1.5\linewidth, 
					]
					\addplot graphics [xmin=0,xmax=768,ymin=0,ymax=819] {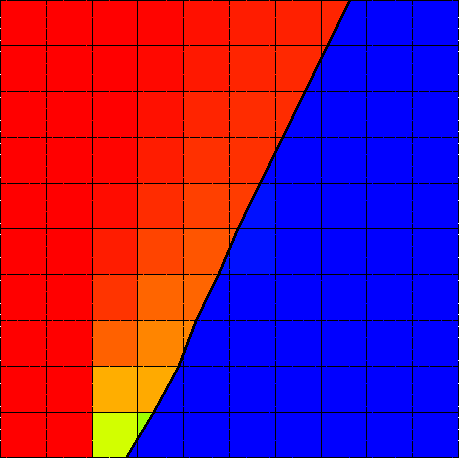};
				\end{axis}
			\end{tikzpicture}
		\end{subfigure}%
		\quad \quad
		\begin{subfigure}{.33\textwidth}
			\centering
			\begin{tikzpicture}
				\begin{axis}[
					title=Iteration 10,
					xlabel={$x$},
					xlabel style={yshift=5pt},
					xmin=0, xmax=768,
					ymin=0, ymax=819,
					xtick={0,384,768},
					xticklabels={0,0.5,1},
					ytick=\empty,
					axis equal image,
					grid=both,
					width=1.5\linewidth, 
					]
					\addplot graphics [xmin=0,xmax=768,ymin=0,ymax=819] {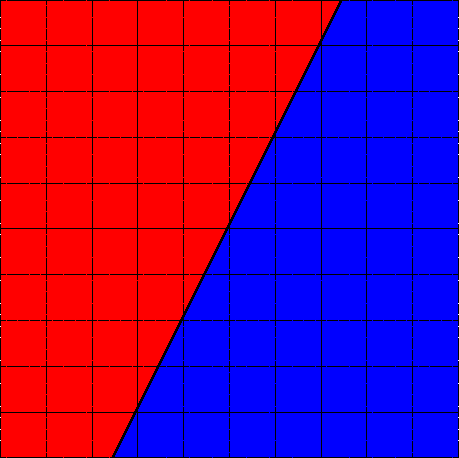};
				\end{axis}
			\end{tikzpicture}
		\end{subfigure}%
		\begin{subfigure}{.33\textwidth}
			\centering
			\begin{tikzpicture}
				\begin{axis}[
					title=Iteration 30,
					xlabel={$x$},
					xlabel style={yshift=5pt},
					xmin=0, xmax=768,
					ymin=0, ymax=819,
					xtick={0,384,768},
					xticklabels={0,0.5,1},
					ytick=\empty,
					axis equal image,
					grid=both,
					width=1.5\linewidth, 
					]
					\addplot graphics [xmin=0,xmax=768,ymin=0,ymax=819] {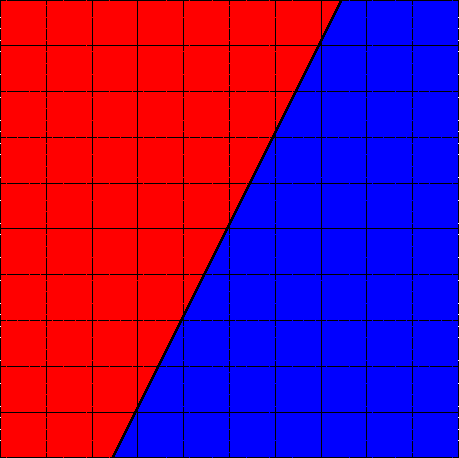};
				\end{axis}
			\end{tikzpicture}
		\end{subfigure}
		\caption{Plots of concentration $c$ for selected XDG shock tracking iterations for the straight shock burgers test case ($P=0$) for $k=0,2,5,10,30$. Here, as an initial guess for the shock interface $\mathfrak{I}_s$ (\textit{thick black line}), implicitly defined by a linear spline level set ($P_s=1$), a curved interface is used which is rapidly aligned to the discontinuity by the SQP solver.}
		\label{fig:SSB}
	\end{figure}
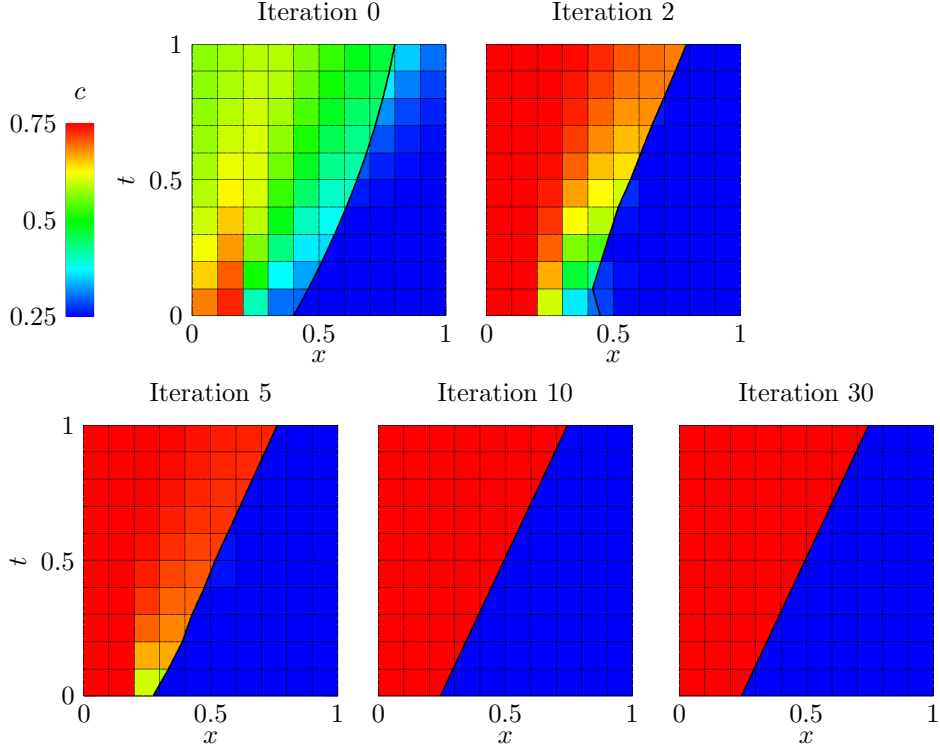
	\begin{figure}[t]
		\centering
		\begin{tikzpicture}
			\begin{semilogyaxis}[            xlabel={Iteration},            ylabel={Residual},            xmin=0, xmax=30,            ymin=1e-16, ymax=10,            xtick={0,10,20,30},            ytick={1,1e-2,1e-4,1e-6,1e-8,1e-10,1e-12,1e-14,1e-16},            xticklabel style={/pgf/number format/.cd,fixed,precision=0},            yticklabel style={/pgf/number format/.cd,fixed,precision=0,sci},            grid=both,            ]
				\addplot table[x index=0,y index=1] {StraightShockBurgers/SSB_EnResNorms.txt};
				\addlegendentry{$\Vert \Rbm(\zbm) \Vert_2$}
				\addplot table[x index=0,y index=1] {StraightShockBurgers/SSB_ResNorms.txt};
				\addlegendentry{$\Vert \rbm(\zbm) \Vert_2$}
			\end{semilogyaxis}
		\end{tikzpicture}
		\caption{Optimization history of the XDG shock tracking for the straight shock burgers test case ($P=0$). Both residual norms $\Vert \Rbm \Vert_2, \Vert\rbm \Vert_2$ converge rapidly after 12 iterations while the SQP solver successfully tracks the straight discontinuity using a linear spline level set.}
		\label{fig:OptHisSSB}
	\end{figure}
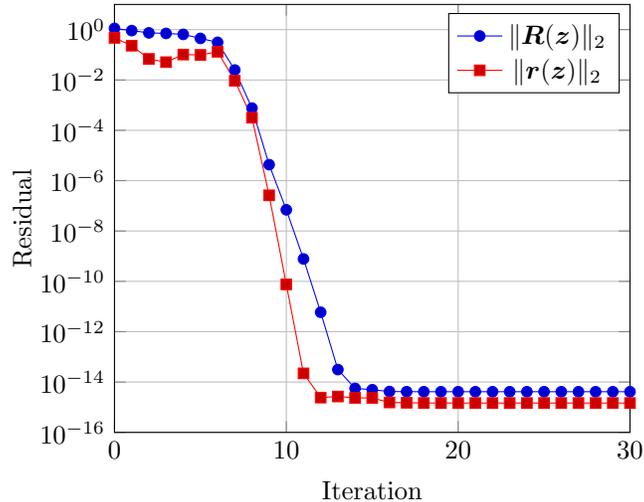 
	Second, the XDG shock tracking method is applied to the non-linear space-time formulation of the inviscid Burgers equation
	\begin{eqnarray}
		\frac{\partial c}{\partial t} + c\frac{\partial c}{\partial x} =0 \quad \text{in }\Omega,
	\end{eqnarray}
	where $ \fbm(c)=(\frac{1}{2}c^2 ~ c)$
	is the corresponding physical flux and $c:\Omega \rightarrow \mathbb{R}$ a concentration. We will investigate two cases: one with a piecewise constant solution and a straight discontinuity, and a second with a curved discontinuity and a non-polynomial solution. For both the upwind flux
	\begin{eqnarray}
		\hat{F}(c^+,c^-,\nbm_\Gamma)= \begin{cases}c^- (\frac{c^++c^-}{2},1)^T \cdot \nbm_\Gamma \quad \text{if } (\frac{c^++c^-}{2},1)^T \cdot \nbm_\Gamma <0 \\
			c^+ (\frac{c^++c^-}{2},1)^T\cdot \nbm_\Gamma) \quad \text{if } (\frac{c^++c^-}{2},1)^T\cdot \nbm_\Gamma \geq 0 \end{cases}.
	\end{eqnarray}
is chosen as the numerical flux.
	\subsubsection{Straight Shock}
	For the first case a rectangular domain $\Omega=[0,1]\times[0,1]$ and a piecewise constant initial value with a discontinuity 
	\begin{eqnarray}
		c(x,0)=\frac{3}{4}  H\left(\frac{1}{4} -x\right) +\frac{1}{4}  H\left(x-\frac{3}{4} \right)
	\end{eqnarray}
are chosen. By applying the Rankine Hugoniot conditions, we determine the shock speed to be $\frac{1}{2}$. This results in a linear shock profile and the exact piecewise constant solution is
	\begin{eqnarray} c_\text{Ex}(x,t)= \frac{3}{4} H\left(-x+\frac{1}{4}+\frac{1}{2}t\right) + \frac{1}{4} H\left(x-\frac{1}{4}-\frac{1}{2}t \right).
	\end{eqnarray}
	Dirichlet boundary conditions are directly prescribed using the exact solution. The solution can be represented in an XDG space of order $P=0$ using the level set
	\begin{eqnarray}
		\varphi_s(x,t)=x-s(t)=x-\frac{1}{4} -\frac{1}{2}t.
	\end{eqnarray} Analogously to the advection test case, a $10 \times 10$ background mesh is chosen. For the straight shock a linear spline level set is used. This example aims to demonstrate the method's capability to solve a non-linear equation and align an initially curved interface with a straight on. Hence, for the initial guess the curved level set
	\begin{eqnarray}
		\varphi_s(x,t)=x-s(t)=x+\frac{1}{5}t^2-\frac{3}{5}t-\frac{2}{5},
	\end{eqnarray} is chosen and projected onto the linear spline. The initial solution guess is obtained from a single Newton step using the fixed initial level set and is illustrated in Figure \ref{fig:SSB}. 

Finally we apply the XDG shock tracking method and after running 30 iterations, we observe convergence at iteration 12. During this process, the discontinuity is successfully tracked and the residuals are reduced to $10^{-14}$. These results demonstrate that the method can effectively solve a non-linear 2D problem with a polynomial order of $P=0$, achieving the correct solution in approximately 12 steps. Selected iterations of the SQP solver are displayed in Figure \ref{fig:SSB}, while Figure \ref{fig:OptHisSAIDT} provides the history of both residual norms.
	\subsubsection{Accelerating shock}
	Next, we demonstrate a case where the solution is no longer piecewise constant and a $P$-continuation strategy is applied. Following Huang \& Zahr \cite{huangRobustHighorderImplicit2022a} the domain $\Omega = [-0.2,1] \times [0,1]$ and the initial condition 
	\begin{equation}
		c(x,0) = \begin{cases}
			4 &\text{if } x<0 \\
			3(x-1) &\text{else}
		\end{cases}
	\end{equation}
	are considered. Here, the analytical solution is known to be 
	\begin{equation}
		c(x,t) = \begin{cases}
			4 &\text{if } x<s(t) \\
			\frac{3(x-1)}{1+3t} &\text{else, }
		\end{cases}
	\end{equation}
	where the accelerating shock speed is given by 
	\begin{equation}
		s(t)= \frac{7}{4}(1-\sqrt{1 +3t}) + 4t.
	\end{equation}
	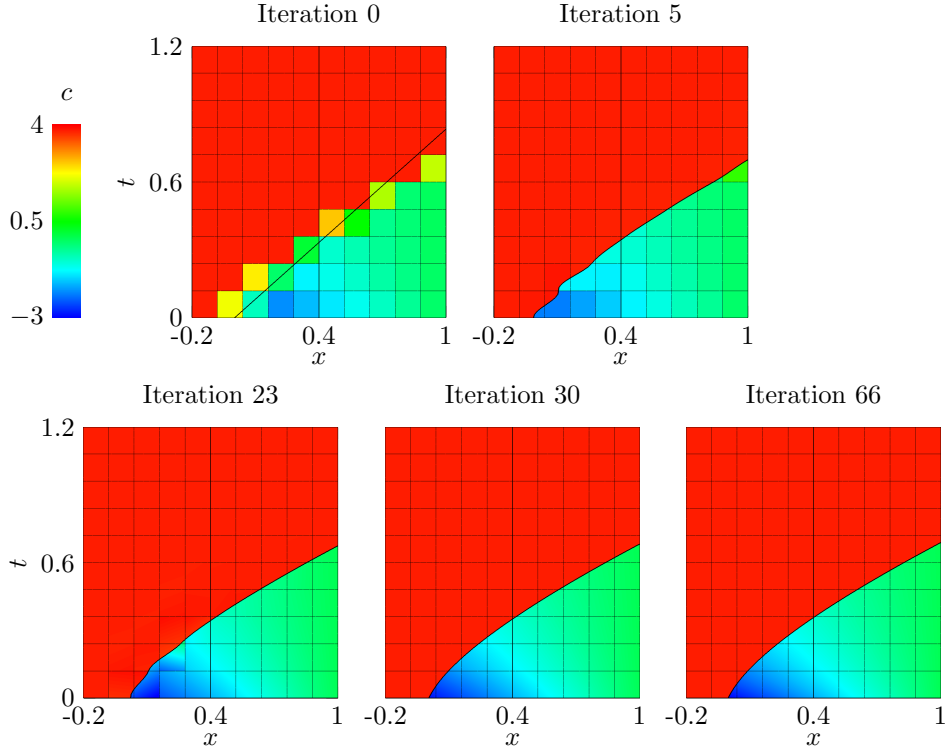
\begin{figure}[t]
		\raisebox{0.5cm}{\begin{subfigure}{0.08\textwidth}
				\centering
				\begin{tikzpicture}
					\begin{axis}[
						title={$c$},
						xlabel={},
						ylabel={},
						xmin=0, xmax=1,
						ymin=-3, ymax=4,
						xtick=\empty,
						ytick={-3,0.5,4},
						grid=both,
						width=4.8cm,  
						axis equal image
						]
						\addplot graphics [xmin=0,xmax=1,ymin=-3,ymax=4] {legend_trimmed.png};
					\end{axis}
				\end{tikzpicture}
		\end{subfigure}}
		\quad
		\begin{subfigure}{.33\textwidth}
			\centering
			\begin{tikzpicture}
				\begin{axis}[
					title=Iteration 0,
					ylabel={$t$},
					xlabel={$x$},
					xlabel style={yshift=5pt}, 
					ylabel style={yshift=-10pt}, 
					xmin=0, xmax=768,
					ymin=0, ymax=819,
					xtick={0,384,768},
					xticklabels={-0.2,0.4,1},
					ytick={0,409,819},
					yticklabels={0,0.6,1.2},
					axis equal image,
					grid=both,
					width=1.5\linewidth, 
					]
					\addplot graphics [xmin=0,xmax=768,ymin=0,ymax=819] {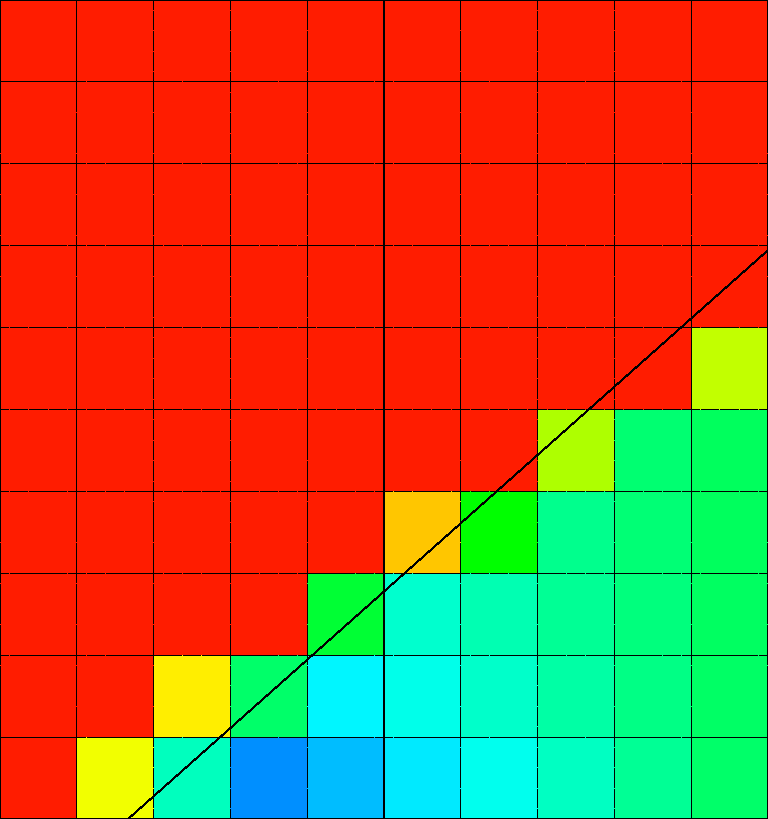};
				\end{axis}
			\end{tikzpicture}
		\end{subfigure}%
		\quad \quad
		\begin{subfigure}{.33\textwidth}
			\centering
			\begin{tikzpicture}
				\begin{axis}[
					title=Iteration 5,
					xlabel={$x$},
					xlabel style={yshift=5pt},
					xmin=0, xmax=768,
					ymin=0, ymax=819,
					xtick={0,384,768},
					xticklabels={-0.2,0.4,1},
					ytick=\empty,
					axis equal image,
					grid=both,
					width=1.5\linewidth, 
					]
					\addplot graphics [xmin=0,xmax=768,ymin=0,ymax=819] {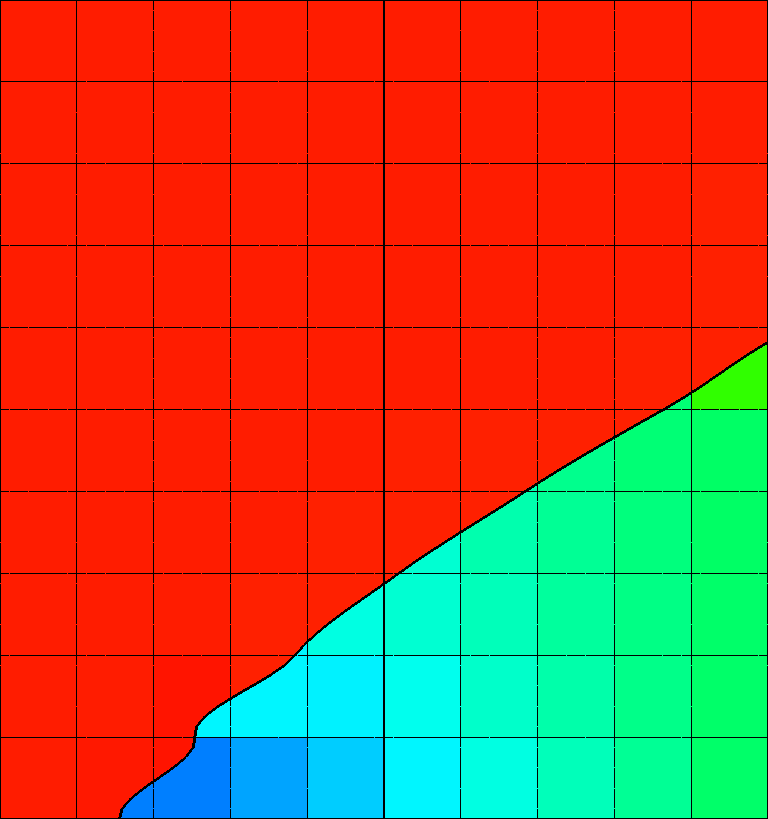};
				\end{axis}
			\end{tikzpicture}
		\end{subfigure}
		\\
		\begin{subfigure}{.33\textwidth}
			\centering
			\begin{tikzpicture}
				\begin{axis}[
					title=Iteration 23,
					xlabel={$x$},
					xlabel style={yshift=5pt},
					xmin=0, xmax=768,
					ymin=0, ymax=819,
					xtick={0,384,768},
					xticklabels={-0.2,0.4,1},
					ytick={0,409,819},
					yticklabels={0,0.6,1.2},
					ylabel={$t$},
					ylabel style={yshift=-10pt},
					axis equal image,
					grid=both,
					width=1.5\linewidth, 
					]
					\addplot graphics [xmin=0,xmax=768,ymin=0,ymax=819] {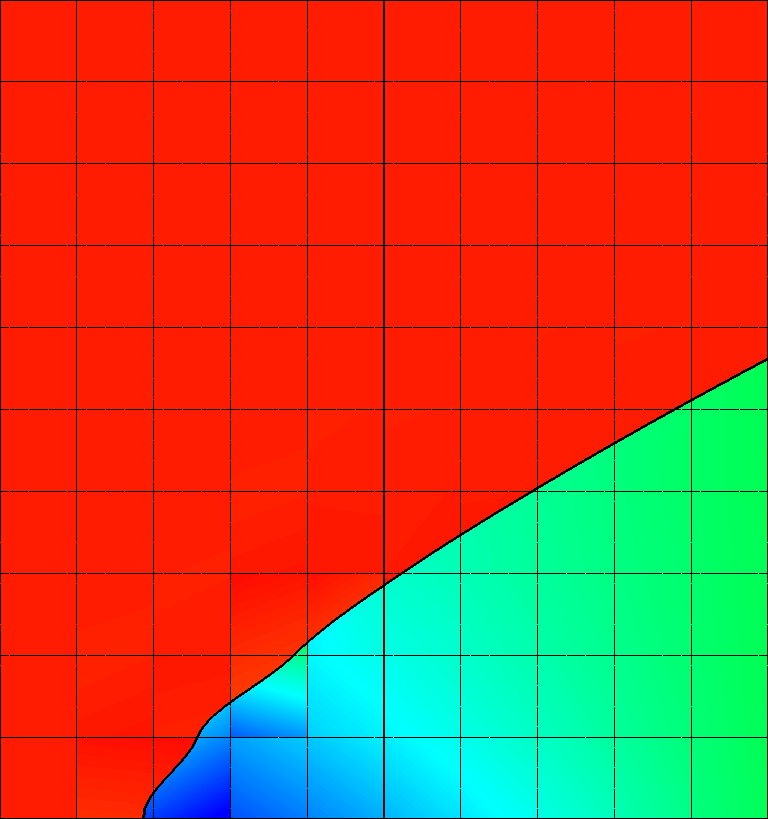};
				\end{axis}
			\end{tikzpicture}
		\end{subfigure}%
			\quad \quad
		\label{fig:ASBfig}
		\begin{subfigure}{.33\textwidth}
			\centering
			\begin{tikzpicture}
				\begin{axis}[
					title=Iteration 30,
					xlabel={$x$},
					xlabel style={yshift=5pt},
					xmin=0, xmax=768,
					ymin=0, ymax=819,
					xtick={0,384,768},
					xticklabels={-0.2,0.4,1},
					ytick=\empty,
					axis equal image,
					grid=both,
					width=1.5\linewidth, 
					]
					\addplot graphics [xmin=0,xmax=768,ymin=0,ymax=819] {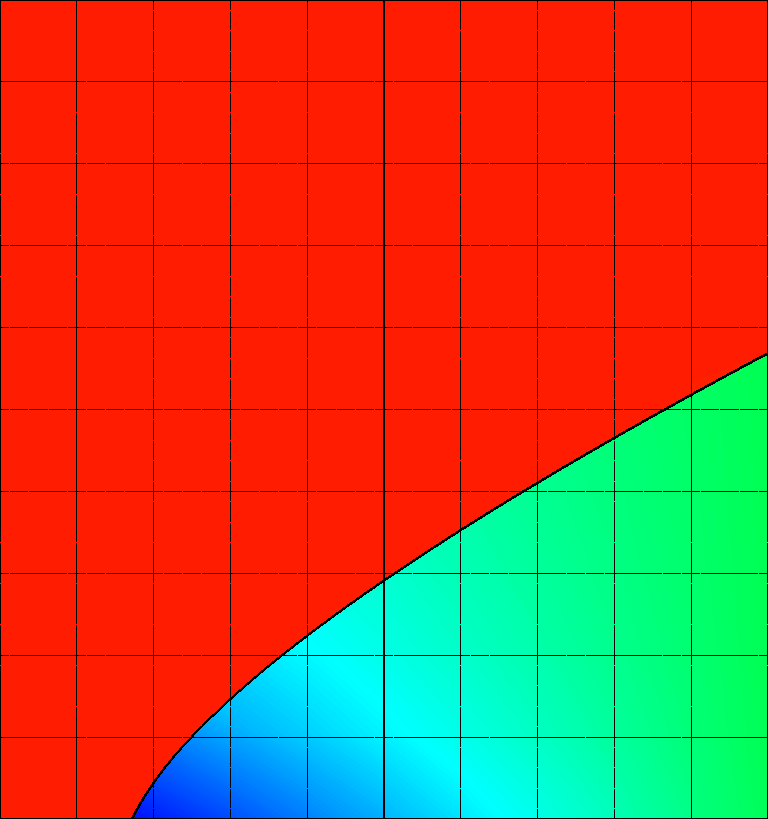};
				\end{axis}
			\end{tikzpicture}
		\end{subfigure}%
		\begin{subfigure}{.33\textwidth}
			\centering
			\begin{tikzpicture}
				\begin{axis}[
					title=Iteration 66,
					xlabel={$x$},
					xlabel style={yshift=5pt},
					xmin=0, xmax=768,
					ymin=0, ymax=819,
					xtick={0,384,768},
					xticklabels={-0.2,0.4,1},
					ytick=\empty,
					axis equal image,
					grid=both,
					width=1.5\linewidth, 
					]
					\addplot graphics [xmin=0,xmax=768,ymin=0,ymax=819] {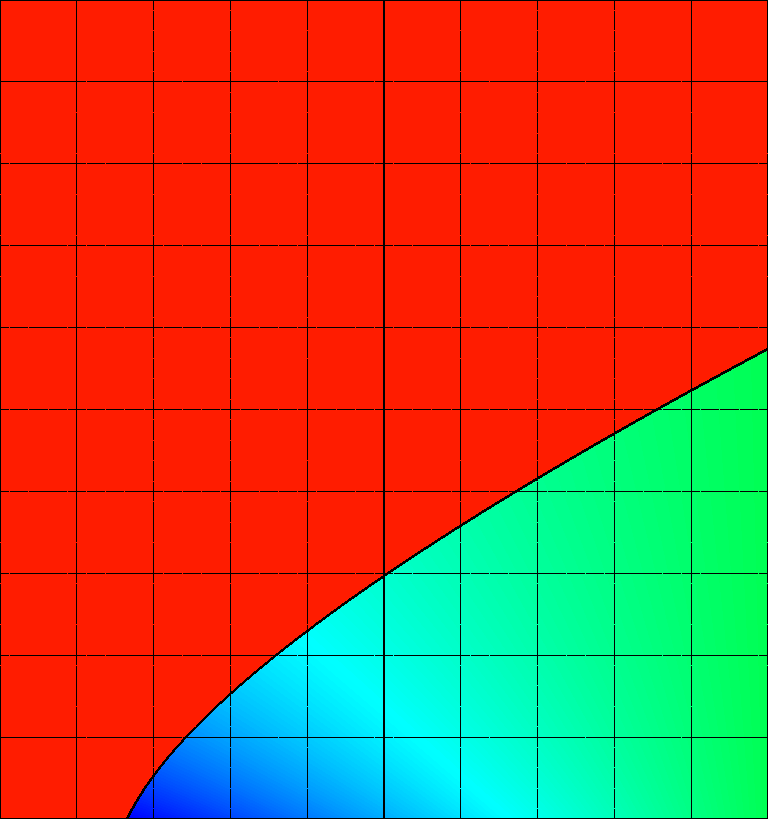};
				\end{axis}
			\end{tikzpicture}
		\end{subfigure}
		\caption{Plots of concentration $c$ for selected XDG shock tracking iterations for the accelerating shock burgers test case for $k=0,5,23,30,66$. Here, as an initial guess for the shock interface $\mathfrak{I}_s$ (\textit{thick black line}), implicitly defined by a cubic spline level set, a linear interface is employed and aligned to the discontinuity by the SQP solver, simultaneously computing the solution of the conservation law. Here, a $P$-continuation strategy is employed gradually increasing the polynomial degrees ($k=1,5 \mapsto P=0;~ k=23 \mapsto P=1;~ k=30 \mapsto P=2;~ k=66 \mapsto P=3$).}
		\label{fig:Burgers2}
	\end{figure}
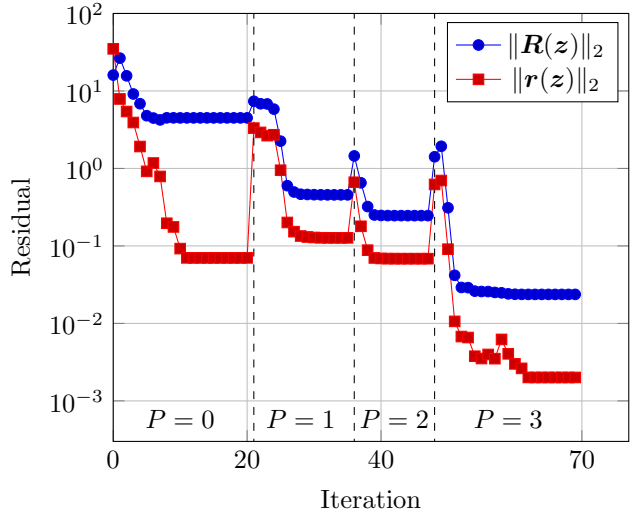
\begin{figure}[t]
	\centering
	\begin{tikzpicture}
		\begin{semilogyaxis}[
			xlabel={Iteration},
			ylabel={Residual},
			xmin=0, xmax=70,
			ymin=3e-4, ymax=100,
			xtick={0,20,40,70},
			ytick={100,10,1,1e-1,1e-2,1e-3,1e-4},
			xticklabel style={/pgf/number format/.cd,fixed,precision=0},
			yticklabel style={/pgf/number format/.cd,fixed,precision=0,sci},
			grid=both,
			enlarge x limits={upper=0.15},
			]
			\addplot table[x index=0,y index=1] {AcceleratingShockBurgers/BurgersIDT_EnResNorms.txt};
			\addlegendentry{$\Vert \Rbm(\zbm) \Vert_2$}
			\addplot table[x index=0,y index=1] {AcceleratingShockBurgers/BurgersIDT_ResNorms.txt};
			\addlegendentry{$\Vert \rbm(\zbm) \Vert_2$}
			
			\draw[dashed] (axis cs:21,1e-4) -- (axis cs:21,100);
			\draw[dashed] (axis cs:36,1e-4) -- (axis cs:36,100);
			\draw[dashed] (axis cs:48,1e-4) -- (axis cs:48,100);
			\node at (axis cs:10,1e-3) [anchor=north] {$P=0$};
			\node at (axis cs:28,1e-3) [anchor=north] {$P=1$};
			\node at (axis cs:42,1e-3) [anchor=north] {$P=2$};
			\node at (axis cs:59,1e-3) [anchor=north] {$P=3$};
		\end{semilogyaxis}
	\end{tikzpicture}
	\caption{Optimization history of the XDG shock tracking for the advection test case where a $P$-continuation strategy is employed ($P=0,1,2,3$) showing the residual norms $\Vert \Rbm \Vert_2, \Vert\rbm \Vert_2$ for all SQP iterations. After converging the residual of an intermediate polynomial degree, the latter is increased, such that a jump in both residuals $\Rbm,\rbm$ can be observed (\textit{dashed lines}). An overall decline of both residuals can be measured as the SQP solver tracks the non-polynomial discontinuity and computes the solution.}
	\label{fig:OptHisASB}
\end{figure}
	This example serves as a test of the methods proficiency in solving a challenging non-linear 2D problem where the solution cannot be precisely represented within the XDG space. Unlike the previous cases, both the interface and the solution are non-polynomial, adding an extra layer of complexity.
	
	Figure \ref{fig:OptHisASB} presents the optimization history, showcasing the $l_2$-norms of the two residuals. As the polynomial degree increases, a noticeable uptick in the residuals can be observed (indicated by dashed lines). After a total of 66 iterations, the method successfully reduces the residual to approximately $3 \times 10^{-3}$, achieving the alignment of the discontinuity in the process. This outcome underscores the method's effectiveness in handling non-linear 2D scenarios with non-polynomial solutions and interfaces.
	\subsection{Inviscid Euler equations}
	\begin{figure}[t]
		\raisebox{0.3cm}{\begin{subfigure}{0.08\textwidth}
				\centering
				\begin{tikzpicture}
					\begin{axis}[
						title={$\rho$},
						xlabel={},
						ylabel={},
						xmin=0, xmax=1,
						ymin=-3, ymax=4,
						xtick=\empty,
						ytick={-3,0.5,4},
						yticklabels={1,1.25,1.5},
						grid=both,
						width=4.8cm,  
						axis equal image
						]
						\addplot graphics [xmin=0,xmax=1,ymin=-3,ymax=4] {legend_trimmed.png};
					\end{axis}
				\end{tikzpicture}
				\label{fig:ASBLegend}
		\end{subfigure}}
		\quad
		\begin{subfigure}{.33\textwidth}
			\centering
			\begin{tikzpicture}
				\begin{axis}[
					title=Iteration 0,
					ylabel={$y$},
					xlabel={$x$},
					xlabel style={yshift=5pt}, 
					ylabel style={yshift=-10pt}, 
					xmin=0, xmax=718,
					ymin=0, ymax=479,
					xtick={0,239,478,718},
					xticklabels={0,0.5,1,1.5},
					ytick={0,239,479},
					yticklabels={0,0.5,1},
					axis equal image,
					grid=both,
					width=1.5\linewidth, 
					]
					\addplot graphics [xmin=0,xmax=718,ymin=0,ymax=479] {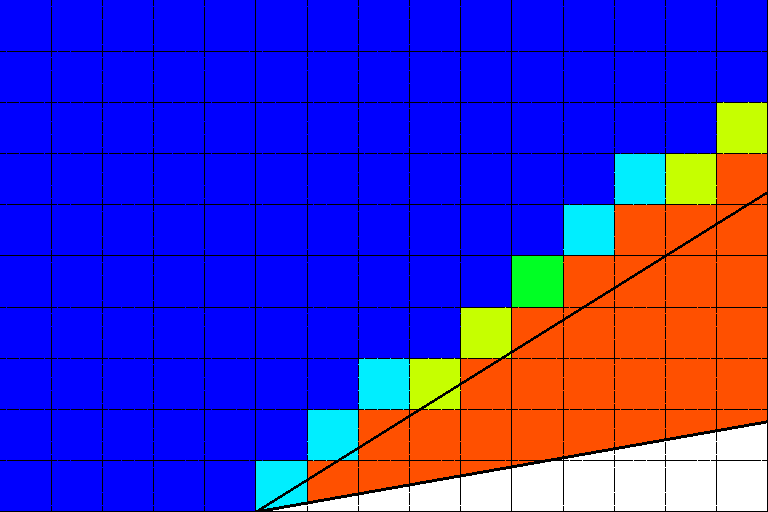};
				\end{axis}
			\end{tikzpicture}
		\end{subfigure}%
		\quad \quad \quad \quad \quad \quad
		\begin{subfigure}{.33\textwidth}
			\centering
			\begin{tikzpicture}
				\begin{axis}[
					title=Iteration 5,
					xlabel={$x$},
					xlabel style={yshift=5pt},
					xmin=0, xmax=718,
					ymin=0, ymax=479,
					xtick={0,239,478,718},
					xticklabels={0,0.5,1,1.5},
					ytick=\empty,
					axis equal image,
					grid=both,
					width=1.5\linewidth, 
					]
					\addplot graphics [xmin=0,xmax=718,ymin=0,ymax=479] {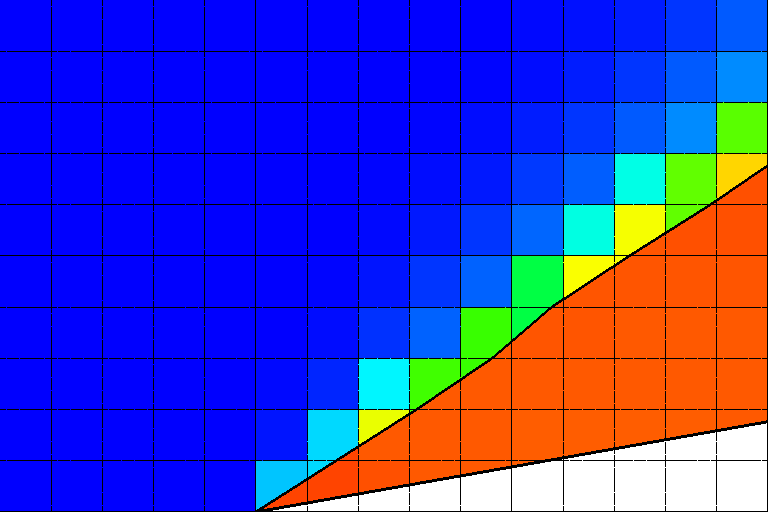};
				\end{axis}
			\end{tikzpicture}
		\end{subfigure}
		\\
		\begin{subfigure}{.33\textwidth}
			\centering
			\begin{tikzpicture}
				\begin{axis}[
					title=Iteration 10,
					ylabel={$y$},
					ylabel style={yshift=-10pt},
					xlabel={$x$},
					xlabel style={yshift=5pt},
					xmin=0, xmax=718,
					ymin=0, ymax=479,
					xtick={0,239,478,718},
					xticklabels={0,0.5,1,1.5},
					ytick={0,239,479},
					yticklabels={0,0.5,1},
					axis equal image,
					grid=both,
					width=1.5\linewidth, 
					]
					\addplot graphics [xmin=0,xmax=718,ymin=0,ymax=479] {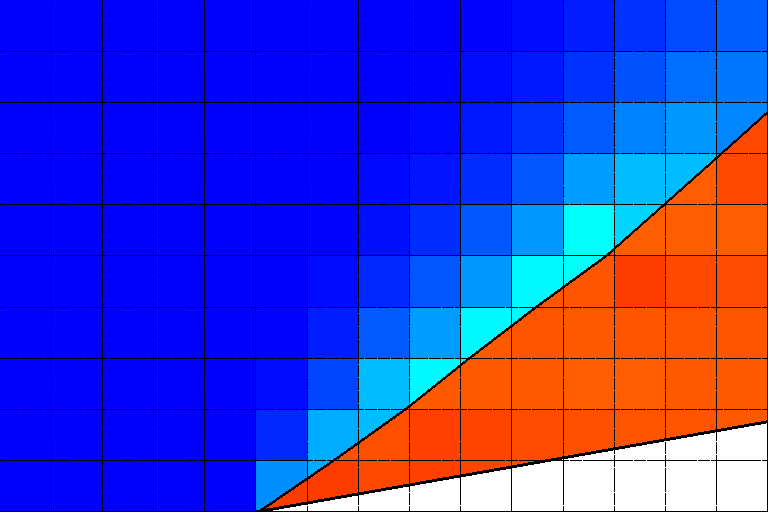};
				\end{axis}
			\end{tikzpicture}
		\end{subfigure}%
		\quad \quad \quad \quad \quad \quad
		\begin{subfigure}{.33\textwidth}
			\centering
			\begin{tikzpicture}
				\begin{axis}[
					title=Iteration 45,
					xlabel={$x$},
					xlabel style={yshift=5pt},
					xmin=0, xmax=718,
					ymin=0, ymax=479,
					xtick={0,239,478,718},
					xticklabels={0,0.5,1,1.5},
					ytick=\empty,
					axis equal image,
					grid=both,
					width=1.5\linewidth, 
					]
					\addplot graphics [xmin=0,xmax=718,ymin=0,ymax=479] {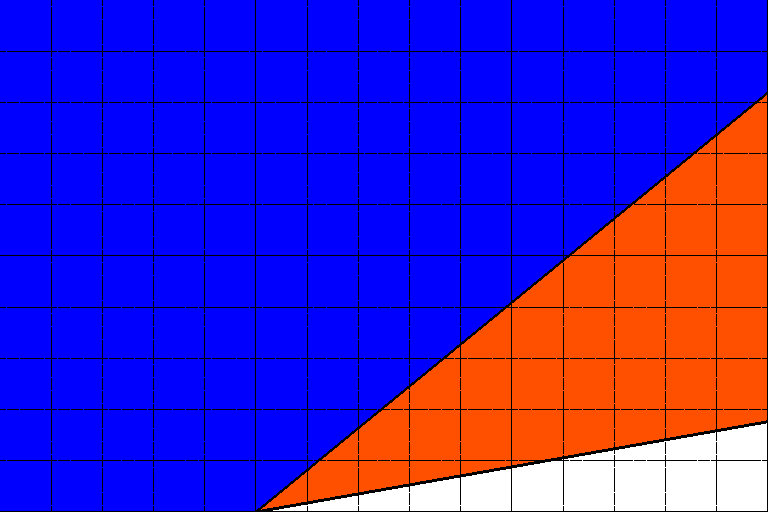};
				\end{axis}
			\end{tikzpicture}
		\end{subfigure}	
		\caption{Selected XDG shock tracking iterations ($k=0,5,10,15$) for the density $\rho$ ($P=0$) of a Mach $2$ flow over an inclined plane (\textit{indicated by white filling}) with angle $\theta_{\text{wedge}}=10^{\circ}$ for a  $15\times10$ mesh, represented by an immersed boundary $\mathfrak{I}_b$ (\textit{lower thick black line},$P_b=1$). The shock interface $\mathfrak{I}_s$ (\textit{upper thick black line}) is represented by a linear spline level set $\varphi_s$ ($P_s=1$) and is converged by the SQP solver to the correct shock position in about 30 iterations.}
		\label{fig:WF2}
	\end{figure}
	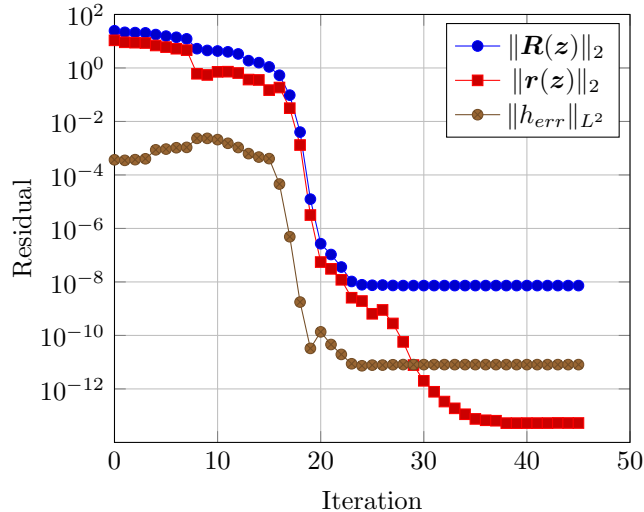
\begin{figure}[t]
		\centering
		\begin{tikzpicture}
			\begin{semilogyaxis}[            xlabel={Iteration},            ylabel={Residual},            xmin=0, xmax=50,            ymin=1e-14, ymax=1e2,            xtick={0,10,20,30,40,50},            ytick={100,1,1e-2,1e-4,1e-6,1e-8,1e-10,1e-12},            xticklabel style={/pgf/number format/.cd,fixed,precision=0},            yticklabel style={/pgf/number format/.cd,fixed,precision=0,sci},            grid=both,            ]
				\addplot table[x index=0,y index=1] {WedgeFlow/R1.txt};
				\addlegendentry{$\Vert \Rbm(\zbm) \Vert_2$}
				\addplot table[x index=0,y index=1] {WedgeFlow/R0.txt};
				\addlegendentry{$\Vert \rbm(\zbm) \Vert_2$}
				\addplot table[x index=0,y index=1] {WedgeFlow/herr.txt};
				\addlegendentry{$\Vert h_{err} \Vert_{L^2}$}
			\end{semilogyaxis}
		\end{tikzpicture}
		\caption{Optimization history of the XDG shock tracking for the Mach $2$ wedge flow ($P=0$) showing the residual norms $\Vert \Rbm \Vert_2, \Vert\rbm \Vert_2$ and the normalized enthalpy error $\Vert h_{\text{err}} \Vert_{L^2}$ across all SQP iterations. An convergence of both residuals and the enthalpy error can be measured after 30 iterations as the SQP solver tracks the straight discontinuity and computes the constant solution.}
		\label{fig:OptHisWF2}
	\end{figure} 
	Finally, we apply the method to two stationary supersonic flow problems. These scenarios involve solving the steady two-dimensional Euler equations governing compressible inviscid flow, which are derived from conservation laws for mass, momentum, and energy \cite{andersonModernCompressibleFlow2003}. The governing equations can be expressed as
	\begin{equation}
		\label{eqn:euler}  \pder{\fbm^x(U)}{x} + \pder{\fbm^y(U)}{y} = 0,
	\end{equation}
	with the coressponding physical flux to be $\fbm(U)=(\fbm^x(U)~ \fbm^y(U))$. Here, the conserved quantities $U: \Omega \rightarrow \mathbb{R}^4$, the convective fluxes $\fbm^x: \mathbb{R}^4 \rightarrow \mathbb{R}^4$ and $\fbm^y: \mathbb{R}^4 \rightarrow \mathbb{R}^4$ are given by the following expressions:
	\begin{equation}
		U = \begin{pmatrix}
			\rho\\
			\rho u\\
			\rho v\\
			\rho E
		\end{pmatrix},\quad
		\fbm^x(U) = \begin{pmatrix}
			\rho u\\
			\rho u^2 + p\\
			\rho u v\\
			u (\rho E + p)
		\end{pmatrix},\quad
		\fbm^y(U) = \begin{pmatrix}
			\rho v\\
			\rho v u\\
			\rho v^2 + p\\
			v (\rho E + p)
		\end{pmatrix},
	\end{equation}
	where $x \in \Omega $, $\rho: \Omega \rightarrow \mathbb{R}^+$ is the fluid density, $\rho u, \rho v: \Omega \rightarrow \mathbb{R}$ are the momentums, $u, v: \Omega \rightarrow \mathbb{R}$ are velocities, $\rho E: \Omega \rightarrow \mathbb{R}^+$ is the total energy per volume and $p: \Omega \rightarrow \mathbb{R}^+$ is the pressure. The total energy itself is the sum of the inner energy $\rho e_{in}: \Omega \rightarrow \mathbb{R}^+$ and the kinetic energy
	\begin{equation}
		\rho E
		= \rho e_{in} + \frac{1}{2} \rho \left(uu+vv\right).
	\end{equation}
	To complete the equations, an equation of state based on the ideal gas law is employed, defined as
	\begin{equation}
		p=(\gamma -1)(\rho E - \frac{1}{2}\rho \left(uu+vv\right)),
	\end{equation}
	where the heat capacity is assumed to be $\gamma = 1.4$ for air at standard conditions.
	As for the numerical flux function $\hat{F}$, the following choices are made for cell boundaries $\Gamma$:
	\begin{itemize}
		\item The Harten-Lax-van Leer-compact (HLLC) flux is used for all boundaries that do not lie on the interface(i.e. $(\Gamma \setminus \mathfrak{I}_s) \setminus \mathfrak{I}_b$) while directly imposing boundary conditions for edges corresponding to the domain boundary. This includes supersonic inlet conditions for the left side of the domain, supersonic outlet conditions for the right side, and adiabatic slip walls for the top and bottom boundaries.
		\item The Godunov flux is utilized for edges corresponding to the shock interface $\mathfrak{I}_{s}$.
		\item An adiabatic slip wall condition is applied for the edges corresponding to the immersed boundary interface $\mathfrak{I}_{b}$.
	\end{itemize} 
	We choose the Godunov flux for the shock interface, because the standard HLLC flux is inconsistent for interfaces with jumps. While the Godunov flux is typically more accurate, it can be computationally expensive. However, given that the fluxes across interfaces constitute a minority of the total flux computations, this choice does not significantly impact performance and is thus justified for the considered cases.
	
	To assess the accuracy of the numerical solutions, the resulting enthalpy is compared to the free-stream enthalpy $h_\infty$, which is known to remain constant for stationary inviscid compressible flows. It can be calculated from prescribed values at the inflow (i.e. $\rho_{\text{in}},p_{\text{in}},M_\infty$) by
	\begin{equation}
	h_\infty= \frac{\gamma}{\gamma-1} + \frac{1}{2}M_\infty^2 \gamma \frac{p_{\text{in}}}{\rho_{\text{in}}}
	\end{equation}
when no inflow in $y$-direction is assumed, i.e. $v_{\text{in}}=0$.
	 Any deviation from this state serves as a metric for the accuracy of the numerical method. Consequently, the normalized $L^2$-error in enthalpy within the fluid domain is computed as
	\begin{equation}
		\Vert h_{\text{err}}\Vert_{L^2} := \frac{1}{\Vert h_\infty \Vert_{L^2}}\sqrt{\int_{\Omega} (h - h_\infty)^2 dV},
	\end{equation} where $h$ denotes the enthalpy computed from the numerical solution.
	
	\subsubsection{Supersonic flow over an inclined plane}
	The first test case presented involves a stationary flow over a wedge, which has been previously used as a test case by Zahr et al. \cite{zahrImplicitShockTracking2020}. It considers a rectangular domain $\Omega = [0,\frac{3}{2}] \times [0,1]$ and a supersonic flow over an inclined plane with an angle of $\theta_{\text{wedge}} >0$. The incoming uniform flow, combined with the geometry of the wedge, results in a piecewise constant flow solution, and the ensuing shock is straight. In fact, in this specific configuration, it is possible to explicitly calculate the angle of the shock wave and the solution using relations presented in  \cite{andersonModernCompressibleFlow2003}. For an inflow condition with Mach number $M_\infty=2$ and a wedge angle of $\theta_{\text{wedge}}=10^{\circ}$, the angle of the shock wave $\theta_{\text{shock}}\approx 39.31^{\circ}$ can be obtained. Furthermore, the inflow pressure $p_\text{in}=1$, inflow density $\rho_\text{in}=1$ and inflow velocities $u_\text{in}=M_\infty\sqrt{\gamma}=2\sqrt{\gamma}$ and $v_\text{in}=0$ are prescribed. 
	
	In this setup, two level sets are utilized: a fixed one, denoted as $\varphi_{b}$, representing the immersed boundary and a variable one, denoted as $\varphi_{s}$, representing the shock. The wedge's geometry can be accurately represented by the zero iso-contour of the following linear level set:
	\begin{equation}
		\varphi_{b}(x,y)=\frac{1}{2} +\frac{y}{\tan(\frac{10^{\circ} \pi}{180^{\circ}})} -x.
	\end{equation}
	Similarly, the exact shock would be approximately represented by the linear level set 
	\begin{equation}
		\varphi_{s}(x,y)=\frac{1}{2} +\frac{y}{\tan(\frac{39.31^{\circ} \pi}{180^{\circ}})} -x.
	\end{equation}
	For the domain discretization a $15 \times 10$ background grid is utilized and an initial projection of the exact solution is employed for the initial guess. Due to the piecewise constant nature of the flow, a polynomial degree of $P=0$ is sufficient for representing the solution, along with linear level sets ($P_s=P_b=1$). Thus, a linear spline level set with an initial angle of $32^{\circ}$ is chosen, intentionally creating a configuration where the level set is not sub-cell accurate. 
	
	Remarkably, applying the XDG shock tracking method demonstrates impressive convergence, with the residual norm reaching approximately $1 \times 10^{-12}$ after about 30 iterations. Furthermore, the enthalpy error decreases to around $10^{-11}$, with the exact enthalpy $h_\infty=6.3$ serving as the reference. Consequently, it can be reasonably inferred that the shock is accurately tracked, and the correct solution is computed with errors primarily attributed to rounding.
	
	 It is worth noting that this particular test case involves doubly cut-cells where the wedge and the shock intersect. These cells pose challenges for integration, and recently developed algorithms by Saye \cite{sayeHighorderQuadratureMulticomponent2022} and Beck et al. \cite{beck_high-order_2023}, the latter of which is implemented in BoSSS, have been instrumental in handling such cases. These algorithms build upon the same quadrature rules that we utilize for single cut-cells \cite{sayeHighOrderQuadratureMethods2015}. 
	\subsubsection{Mach 4 bow shock}
	\begin{figure}[t]
		\raisebox{2.4cm}{\begin{subfigure}{0.08\textwidth}
				\centering
				\begin{tikzpicture}
					\begin{axis}[
						title={$\rho$},
						xlabel={},
						ylabel={},
						xmin=0, xmax=1,
						ymin=-3, ymax=4,
						xtick=\empty,
						ytick={-3,0.5,4},
						yticklabels={1,3.05,5.1},
						grid=both,
						width=4.8cm,  
						axis equal image
						]
						\addplot graphics [xmin=0,xmax=1,ymin=-3,ymax=4] {legend_trimmed.png};
					\end{axis}
				\end{tikzpicture}
				\label{fig:BowShockLegend}
		\end{subfigure}}
		\quad
		\begin{subfigure}{.21\textwidth}
			\centering
			\begin{tikzpicture}
				\begin{axis}[
					title=Iteration 0,
					ylabel={$y$},
					xlabel={$x$},
					xlabel style={yshift=5pt}, 
					ylabel style={yshift=-20pt}, 
					xmin=0, xmax=206,
					ymin=0, ymax=822,
					xtick={2,103,206},
					xticklabels={-2,-1,0},
					ytick={0,205,411,617,822},
					yticklabels={-4,-2,0,2,4},
					axis equal image,
					grid=both,
					width=4\linewidth, 
					]
					\addplot graphics [xmin=0,xmax=206,ymin=0,ymax=822] {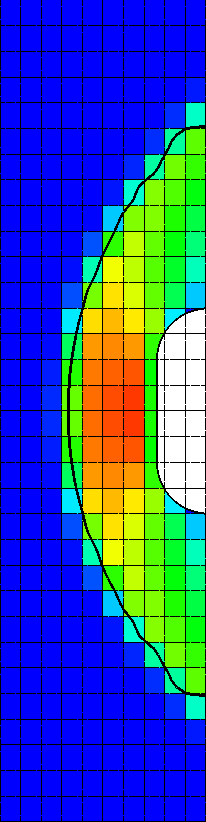};
				\end{axis}
			\end{tikzpicture}
		\end{subfigure}%
		\quad
		\begin{subfigure}{.21\textwidth}
			\centering
			\begin{tikzpicture}
				\begin{axis}[
					title=Iteration 4,
					ylabel={},
					xlabel={$x$},
					xlabel style={yshift=5pt}, 
					xmin=0, xmax=206,
					ymin=0, ymax=822,
					xtick={2,103,206},
					xticklabels={-2,-1,0},
					ytick=\empty,
					axis equal image,
					grid=both,
					width=4\linewidth, 
					]
					\addplot graphics [xmin=0,xmax=206,ymin=0,ymax=822] {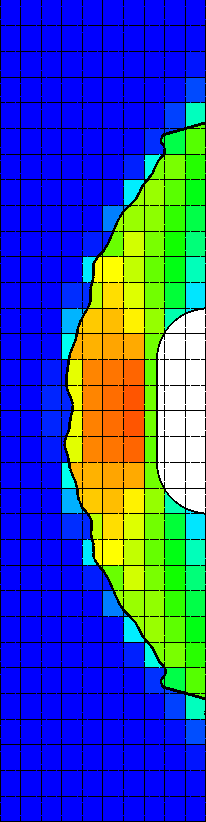};
				\end{axis}
			\end{tikzpicture}
		\end{subfigure}%
		\begin{subfigure}{.21\textwidth}
			\centering
			\begin{tikzpicture}
				\begin{axis}[
					title=Iteration 8,
					ylabel={},
					xlabel={$x$},
					xlabel style={yshift=5pt}, 
					xmin=0, xmax=206,
					ymin=0, ymax=822,
					xtick={2,103,206},
					xticklabels={-2,-1,0},
					ytick=\empty,
					axis equal image,
					grid=both,
					width=4\linewidth, 
					]
					\addplot graphics [xmin=0,xmax=206,ymin=0,ymax=822] {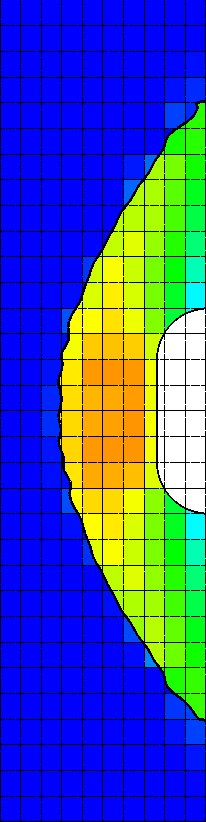};
				\end{axis}
			\end{tikzpicture}
		\end{subfigure}%
		\begin{subfigure}{.21\textwidth}
			\centering
			\begin{tikzpicture}
				\begin{axis}[
					title=Iteration 34,
					ylabel={},
					xlabel={$x$},
					xlabel style={yshift=5pt}, 
					xmin=0, xmax=206,
					ymin=0, ymax=822,
					xtick={2,103,206},
					xticklabels={-2,-1,0},
					ytick=\empty,
					axis equal image,
					grid=both,
					width=4\linewidth, 
					]
					\addplot graphics [xmin=0,xmax=206,ymin=0,ymax=822] {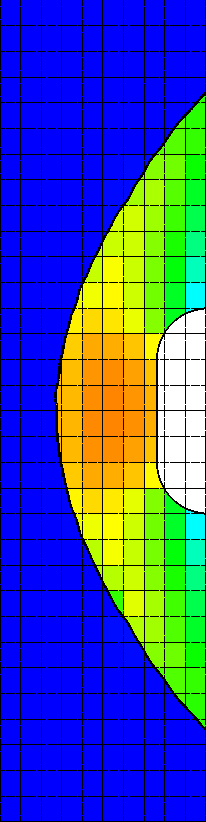};
				\end{axis}
			\end{tikzpicture}
		\end{subfigure}%
		\caption{Selected XDG shock tracking iterations ($P=0$) for the density $\rho$ of the Mach $4$ bow shock for the $10\times32$ mesh and for $k=0,4,8,34$. The blunt body (\textit{white filling}) is represented by an immersed boundary $\mathfrak{I}_b$ (\textit{right thick black line}) using a quadratic level set. Starting from an initial guess obtained by a reconstruction procedure, the shock interface $\mathfrak{I_s}$ (\textit{left thick black line}), represented by a cubic spline level set ($P_s=3$), and the solution are converged by the SQP solver.}
		\label{fig:BowShockP0}
	\end{figure}
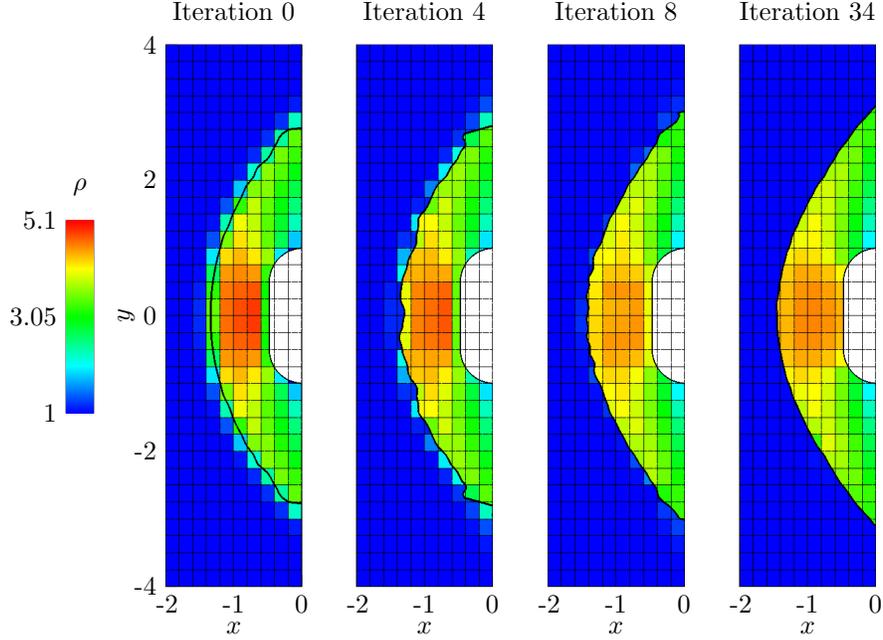
	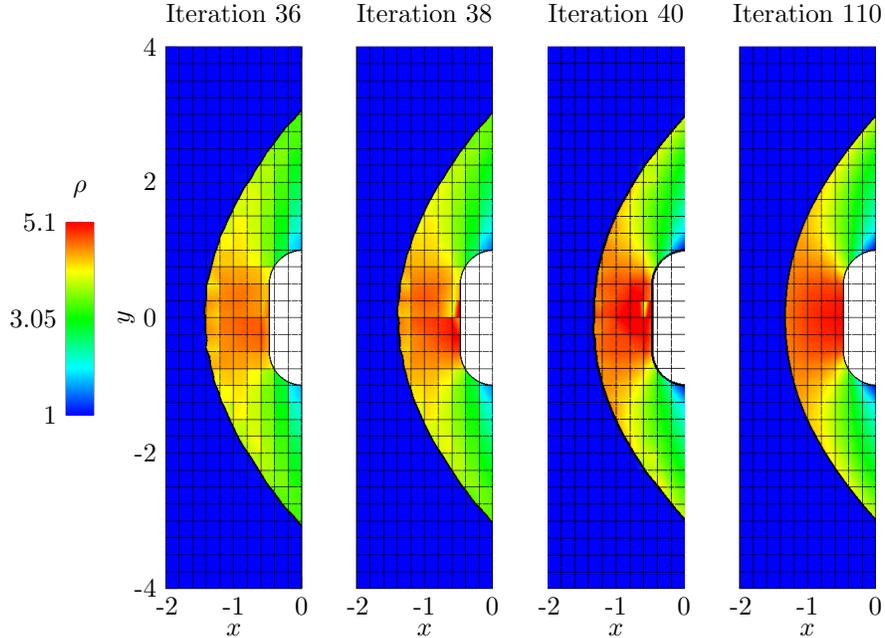
\begin{figure}[t]
		\raisebox{2.4cm}{\begin{subfigure}{0.08\textwidth}
				\centering
				\begin{tikzpicture}
					\begin{axis}[
												title={$\rho$},
						xlabel={},
						ylabel={},
						xmin=0, xmax=1,
						ymin=-3, ymax=4,
						xtick=\empty,
						ytick={-3,0.5,4},
						yticklabels={1,3.05,5.1},
						grid=both,
						width=4.8cm,  
						axis equal image
						]
						\addplot graphics [xmin=0,xmax=1,ymin=-3,ymax=4] {legend_trimmed.png};
					\end{axis}
				\end{tikzpicture}
				\label{fig:BowShock2Legend}
		\end{subfigure}}
		\quad
		\begin{subfigure}{.21\textwidth}
			\centering
			\begin{tikzpicture}
				\begin{axis}[
					title=Iteration 36,
					ylabel={$y$},
					xlabel={$x$},
					xlabel style={yshift=5pt}, 
					ylabel style={yshift=-20pt}, 
					xmin=0, xmax=206,
					ymin=0, ymax=822,
					xtick={2,103,206},
					xticklabels={-2,-1,0},
					ytick={0,205,411,617,822},
					yticklabels={-4,-2,0,2,4},
					axis equal image,
					grid=both,
					width=4\linewidth, 
					]
					\addplot graphics [xmin=0,xmax=206,ymin=0,ymax=822] {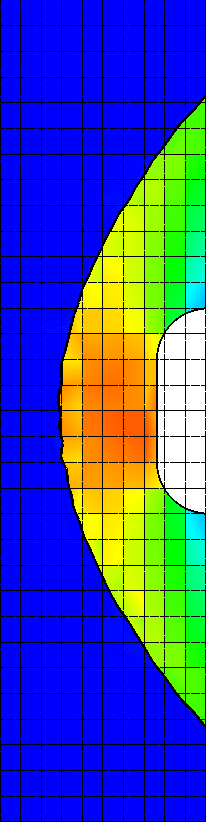};
				\end{axis}
			\end{tikzpicture}
		\end{subfigure}%
		\quad
		\begin{subfigure}{.21\textwidth}
			\centering
			\begin{tikzpicture}
				\begin{axis}[
					title=Iteration 38,
					ylabel={},
					xlabel={$x$},
					xlabel style={yshift=5pt}, 
					xmin=0, xmax=206,
					ymin=0, ymax=822,
					xtick={2,103,206},
					xticklabels={-2,-1,0},
					ytick=\empty,
					axis equal image,
					grid=both,
					width=4\linewidth, 
					]
					\addplot graphics [xmin=0,xmax=206,ymin=0,ymax=822] {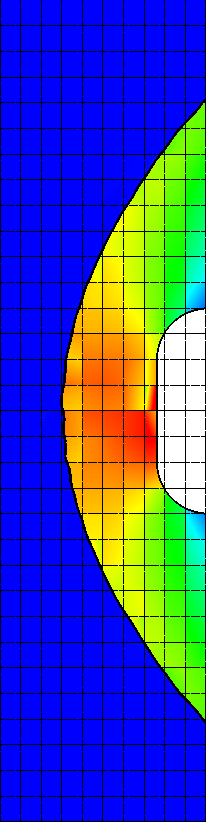};
				\end{axis}
			\end{tikzpicture}
		\end{subfigure}%
		\begin{subfigure}{.21\textwidth}
			\centering
			\begin{tikzpicture}
				\begin{axis}[
					title=Iteration 40,
					ylabel={},
					xlabel={$x$},
					xlabel style={yshift=5pt}, 
					xmin=0, xmax=206,
					ymin=0, ymax=822,
					xtick={2,103,206},
					xticklabels={-2,-1,0},
					ytick=\empty,
					axis equal image,
					grid=both,
					width=4\linewidth, 
					]
					\addplot graphics [xmin=0,xmax=206,ymin=0,ymax=822] {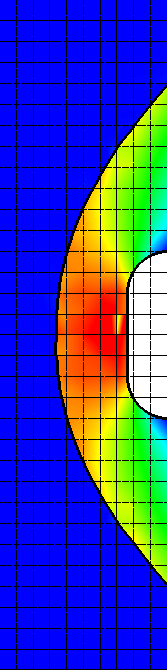};
				\end{axis}
			\end{tikzpicture}
		\end{subfigure}%
		\begin{subfigure}{.21\textwidth}
			\centering
			\begin{tikzpicture}
				\begin{axis}[
					title=Iteration 110,
					ylabel={},
					xlabel={$x$},
					xlabel style={yshift=5pt}, 
					xmin=0, xmax=206,
					ymin=0, ymax=822,
					xtick={2,103,206},
					xticklabels={-2,-1,0},
					ytick=\empty,
					axis equal image,
					grid=both,
					width=4\linewidth, 
					]
					\addplot graphics [xmin=0,xmax=206,ymin=0,ymax=822] {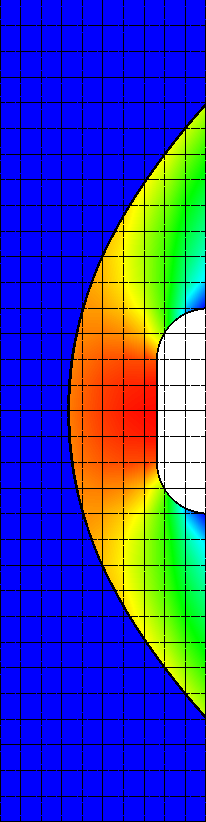};
				\end{axis}
			\end{tikzpicture}
		\end{subfigure}%
		\caption{Selected XDG shock tracking iterations ($P=1$) for the density $\rho$ of the Mach $4$ bow shock for the $10\times32$ mesh and for $k=36,38,40,110$. The blunt body (\textit{white filling}) is represented by an immersed boundary $\mathfrak{I}_b$ (\textit{right thick black line}) using a quadratic level set. Starting from the converged $P=0$ shock tracking solution (see Figure \ref{fig:BowShockP0}), the shock interface $\mathfrak{I_s}$ (\textit{left thick black line}), represented by a cubic spline level set ($P_s=3$), and the solution are converged by the SQP solver.}
		\label{fig:BowShockP1}
	\end{figure}
	
	The final case under consideration involves a blunt body problem characterized by a curved geometry and a curved shock based on the configuration proposed by the HiOCFD5-Workshop \cite{murmanCI1InviscidBow2017}. In this setup, the blunt body is represented by a rounded rectangle formed by the quadrants of two circles with central coordinates $\cbm_1=(x_1,y_1)=(0,0.5)$ and $\cbm_2=(x_2,y_2)=(0,-0.5)^T$, each with a radius of $r=0.5$. Additionally, a vertical line segment is included at $\{x=-0.5, y\in [-0.5,0.5]\}$ to complete the representation of the blunt body. To describe this geometry, a continuous level set function $\varphi_{b} \in C^1(\Omega)$ is used in its quadratic formulation:
	\begin{equation}
		\varphi_{b} (\xbm):=\begin{cases}
			\Vert \xbm - \cbm_1 \Vert_2^2 - 0.5^2 &\text{if} ~ y \geq 0.5 \\
			x^2 - 0.5^2  &\text{if} ~ 0.5 > y \geq -0.5 \\
			\Vert \xbm - \cbm_2 \Vert_2^2 - 0.5^2 &\text{else}.
		\end{cases}
	\end{equation}
	This level set can be precisely represented in the DG space of order two $(P_b=2)$ when the mesh is aligned with the horizontal lines ${y=0.5}$ and ${y=-0.5}$. We prescribe supersonic free-stream conditions with a Mach number of $M_\infty=4.0$, the inflow pressure $p_\text{in}=1$, inflow density $\rho_\text{in}=1$, inflow velocities $u_\text{in}=M_\infty\sqrt{\gamma}=4\sqrt{\gamma}$ and $v_\text{in}=0$. Adiabatic slip-wall boundary conditions are applied at the blunt body.
	
	The computational domain is set to $\Omega =[-2+\epsilon,-\epsilon]\times [-4,4]$, with a slight displacement of the proposed domain in the negative $x$-direction, given by $\epsilon=0.0025$. This adjustment ensures a suitable cut configuration where the interface does not align with the straight vertical portion of the blunt body at $x=-0.5$. To represent the shock interface $\mathfrak{I}_s$, we utilize a cubic spline level set ($P_s=3$) as defined in \eqref{equation:spline}. \begin{figure}[t]
		\centering
		\begin{tikzpicture}
			\begin{semilogyaxis}[            xlabel={Iteration},            ylabel={Residual},            xmin=0, xmax=200,            ymin=1e-6, ymax=20000,            xtick={0,40,80,120,160,200},            ytick={1000,100,1,1e-2,1e-4,1e-6},            xticklabel style={/pgf/number format/.cd,fixed,precision=0},            yticklabel style={/pgf/number format/.cd,fixed,precision=0,sci},            grid=both,            ]
				\addplot+[mark options={mark size=0.5pt,solid, line width=0.5pt}] table[x index=0,y index=1] {BowShock/BowShockEnthalpyResPlot.txt};
				\addlegendentry{$\Vert \Rbm(\zbm) \Vert_2$}
				\addplot+[mark options={mark size=0.5pt,solid, line width=0.5pt}] table[x index=0,y index=2] {BowShock/BowShockEnthalpyResPlot.txt};
				\addlegendentry{$\Vert \rbm(\zbm) \Vert_2$}
				\addplot+[mark options={mark size=0.5pt,solid, line width=0.5pt}] table[x index=0,y index=3] {BowShock/BowShockEnthalpyResPlot.txt};
				\addlegendentry{$\Vert h_{err} \Vert_{L^2}$}
				
				\draw[dashed] (axis cs:35,1e-6) -- (axis cs:35,20000);
				\draw[dashed] (axis cs:111,1e-6) -- (axis cs:111,20000);
				\draw[dashed] (axis cs:144,1e-6) -- (axis cs:144,20000);
				\node at (axis cs:17,1e-5) [anchor=north] {$P=0$};
				\node at (axis cs:72,1e-5) [anchor=north] {$P=1$};
				\node at (axis cs:125,1e-5) [anchor=north] {$P=2$};
				\node at (axis cs:175,1e-5) [anchor=north] {$P=3$};
				
			\end{semilogyaxis}
		\end{tikzpicture}
		\caption{Optimization history of the XDG shock tracking for the Mach $4$ bow shock employing a $P$-continuation strategy is employed ($P=0,1,2,3$). The residual norms $\Vert \Rbm \Vert_2, \Vert\rbm \Vert_2$ and the normalized enthalpy error $\Vert h_{\text{err}} \Vert_{L^2}$ across all SQP iterations are shown. After converging the residual of a intermediate polynomial degree, the latter is increased, such that a jump in both residuals $\Rbm,\rbm$ can be observed (\textit{dashed lines}). An overall decline of both residuals and the enthalpy error can be measured as the SQP solver tracks the non-polynomial discontinuity and solution.}
		\label{fig:OptHisBS}
	\end{figure}
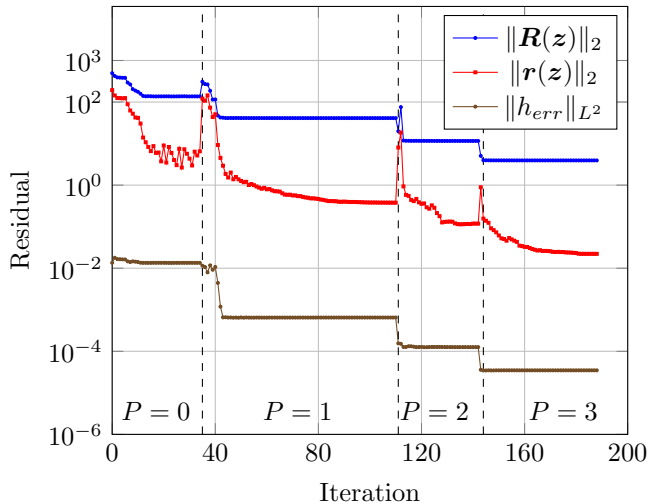
	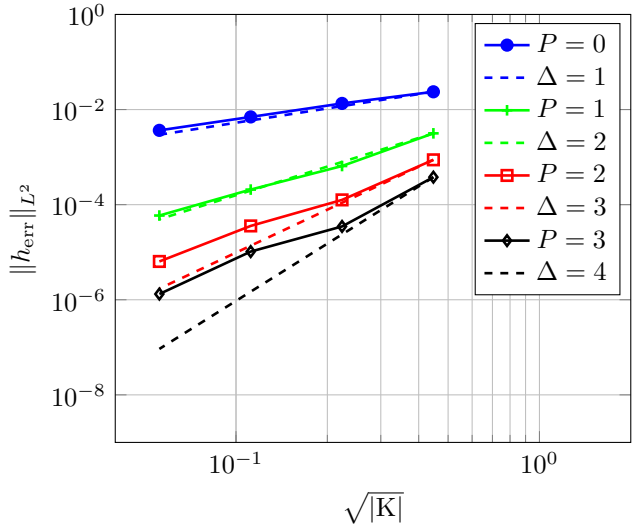
\begin{figure}[t]
\centering
\begin{tikzpicture}
	\begin{loglogaxis}[ xlabel={$\sqrt{\vert \text{K} \vert}$}, ylabel={$ \Vert h_{\text{err}} \Vert_{L^2}$}, xmin=0.04, xmax=2, ymin=1e-09, ymax=1, xtick={1e-1,1}, ytick={1,1e-2,1e-4,1e-6,1e-8}, xticklabel style={/pgf/number format/.cd,fixed,precision=0}, yticklabel style={/pgf/number format/.cd,fixed,precision=0,sci}, grid=both, ]
		\addplot[color=blue, line width=1pt,mark=*] table[x index=0,y index=1] {BowShock/ConvergencePlot.txt};
		\addlegendentry{$P=0$}
		\addplot[color=blue, line width=1pt,dashed] table[x index=0,y index=1] {BowShock/ConvergencePlotIdeal.txt};
		\addlegendentry{$\Delta=1$}
		\addplot[color=green, line width=1pt,mark=+] table[x index=0,y index=2] {BowShock/ConvergencePlot.txt};
		\addlegendentry{$P=1$}
		\addplot[color=green, line width=1pt,dashed] table[x index=0,y index=2] {BowShock/ConvergencePlotIdeal.txt};
		\addlegendentry{$\Delta=2$}
		\addplot[color=red, line width=1pt,mark=square] table[x index=0,y index=3] {BowShock/ConvergencePlot.txt};
		\addlegendentry{$P=2$}
		\addplot[color=red, line width=1pt,dashed] table[x index=0,y index=3] {BowShock/ConvergencePlotIdeal.txt};
		\addlegendentry{$\Delta = 3$}
		\addplot[color=black, line width=1pt,mark=diamond] table[x index=0,y index=4] {BowShock/ConvergencePlot.txt};
		\addlegendentry{$P=3$}
		\addplot[color=black, line width=1pt,dashed] table[x index=0,y index=4] {BowShock/ConvergencePlotIdeal.txt};
		\addlegendentry{$\Delta=4$}
	\end{loglogaxis}
\end{tikzpicture}
\caption{Convergence plot for normalized enthalpy error $\Vert h_{\text{err}} \Vert_{L^2}$ for different polynomial degrees ($P=0,1,2,3$) for the Mach $4$ bow shock on a fixed domain resolved with different mesh resolutions $\sqrt{\vert K \vert}$. Here, the dashed lines show the expected convergence ($P+1$) of a high order method. While $P=0,1$ are in line with expectation $P=2,3$ show a deviation from the expected errors.}
\label{fig:ConvPlot}
\end{figure} 

The initial guess for the numerical solution and the level set are derived from a DG solution of order $P=2$, obtained using explicit time stepping ($t_{\text{end}}=16.0$) and shock capturing on a fine grid composed of $40 \times 160$ cells \cite{geisenhofermarkusShockCapturingHighOrderShockFitting2021}. Employing a projection, the DG solution is used to initialize the XDG solution and a special reconstruction procedure is applied to obtain the initial level set approximation. This procedure aims to find a set of points $(\chi_s:=\{\xbm_i\})$ that describe the entire shock front and can be summarized in the following steps:
\begin{enumerate}
\item Points $\chi_{\text{seed}}$ are seeded into cells where artificial viscosity is applied.
\item From these seed points, candidate points $\chi_{\text{cand}}$, featuring a sign change of the Hessian $H_\rho$ along the density gradient $\nabla \rho$, are computed using a bisection algorithm.
\item Two clustering steps are performed, grouping the points based on their density values and amount of artificial viscosity applied, discarding all points not on the shock.
\end{enumerate}
Further details regarding this procedure can be found in \cite {geisenhofermarkusShockCapturingHighOrderShockFitting2021} and the resulting points are openly accessible \cite{geisenhoferBowShockInflectionpoints2021}.
 These points are then interpolated using the spline level set. Additionally, the $P$-continuation strategy is employed, commencing with $P=0$ and progressively increasing the degree of the solution when stagnation is observed.

We conduct a spatial convergence study for various polynomial degrees $P = {0, 1, 2, 3}$ and discretize the computational domain $\Omega$ using equidistant Cartesian grids comprising $5 \times 16$, $10 \times 32$, $20 \times 64$, and $40 \times 128$ cells. The cell-agglomeration technique is applied with an agglomeration threshold of $\delta_\mathrm{agg} = 0.4$ at both interfaces $\mathfrak{I}_b$ and $\mathfrak{I}_s$. For illustrative purposes, we provide plots for selected iterations for $P=0$ (Figure \ref{fig:BowShockP0}) and for $P=1$ (Figure \ref{fig:BowShockP1}), along with the complete optimization history (Figure \ref{fig:OptHisBS}) for the 10 x 32 mesh.

As described earlier, we assess the accuracy of the numerical solutions by comparing them to the free-stream enthalpy, denoted as $h_\infty$, since an analytical solution for the flow field is unavailable. For this particular test case, $h_\infty = 14.7$ holds consistently throughout the entire flow field. Hence, any deviation from this state serves as a metric for the accuracy of the numerical method. Figure \ref{fig:ConvPlot} presents the results of the spatial convergence study, highlighting the convergence of the methods across various polynomial degrees.
\begin{table}[htbp]
\centering
\caption{Estimated orders of convergence for different polynomial degrees P for the Mach $4$ bow shock on a fixed domain resolved with different mesh resolutions (see Figure \ref{fig:ConvPlot}) }
\label{table:EOC}
\begin{tabular}{ccccc}
	\hline
	$P$ & EOC \\ 
	\hline
	0 & 0.90\\ 
	1 & 1.92 \\
	2 & 2.37\\
	3 & 2.69 \\
	\hline
\end{tabular}
\end{table}
While the estimated rate of convergence, depicted in Table \ref{table:EOC}, 
increases with the polynomial degrees, it is not optimal as it is desired for a high-order method. During this work, a specific reason for this could not be identified but the results show a benefit from using the high-order approximation. 

Overall, we observed that this test case posed significantly greater challenges compared to the previous examples. Even slight modifications to the initial guess of the level set, the state, or variations in the agglomeration parameter could lead to dramatic impacts on the optimization history, particularly noticeable for the coarser meshes ($5\times 16$ and $10\times32$ cells). This hints that the robustness of the method can be greatly improved by an suitable setup procedure which will be the issue of further works. One contributing factor appears to be the behavior of the spline that defines the shock interface, where we observed oscillatory behavior in the sub-iterations. Choosing an alternative level set representation, such as a spline that minimizes curvature, or introducing a penalty term to the objective function, as suggested by Huang \& Zahr \cite{huangRobustHighorderImplicit2022a}, could potentially enhance robustness but were not explored in this study. Another factor may stem from the $P$-continuation strategy. Comparing the converged $P=0$ solution, as depicted in Figure \ref{fig:BowShockP0}, to the converged $P=1$ solution in Figure \ref{fig:BowShockP1}, reveals that the shock interfaces reside in different cells. This may significantly impact convergence as through our simulations we specifically observed that configurations where the interface needs to cross cells can pose challenges for the SQP solver. 
\section{Conclusion \& Outlook}
In this work, we introduced a novel high-order implicit XDG shock tracking method, demonstrating the potential of extending a DG-IBM, using a level set $\varphi_b$ for geometry immersion \cite{geisenhoferDiscontinuousGalerkinImmersed2019}, with an additional level set $\varphi_s$, which is fitted to the shock front. By doing so, the approximation is space is enriched, allowing for an accurate representation of solution discontinuities within cut-cells. This technique has found successful applications in problems governed by conservation laws featuring solution discontinuities.

To obtain high-order shock-aligned XDG solutions, the presented shock tracking method solves a constrained optimization problem using a quasi-Newton SQP solver. Solution extrapolation in newborn cut-cells, adaptive regularization and globalization through line search subroutines are used to improve the stability of the method. The shock level set is constructed from a spline describing the zero iso-contour explicitly, saving computational time and guaranteeing $C^1$-smoothness of the interface. Furthermore, robustness measures were accordingly presented: cell agglomeration, solution re-initialization, and the $P$-continuation strategy. These measures help to further stabilize the method on its way to being used for a wide range of applications.

As a proof of concept, numerical results for two different space-time problems and the steady Euler equations in two dimensions were shown. Our results showcased the method's ability to accurately track both straight-sided and curved discontinuities, even when initiated with non-sub-cell accurate initial guesses. We applied it to supersonic flow over an inclined plane and a blunt body, achieving promising results. Furthermore, a convergence study demonstrated high-order convergence for the supersonic bow shock problem. While our examples were relatively straightforward, featuring single shocks without complex patterns, our framework holds promise for handling more intricate discontinuity scenarios in PDEs, which will be explored in future work.

To enhance the method's versatility and make it suitable for a wide range of high-order shock tracking applications, further research is necessary. The employed level set representation, optimized for computational efficiency and $C^1$ continuity, currently limits us to shocks representable as height functions. Investigating alternative level set representations with greater flexibility and managing associated costs for Jacobians will be a priority. Potential approaches include shock descriptions using B-spline curves/surfaces together with explicit quadrature rules for cut-cells, and cell-local continuous level sets with displacement fields, offering the possibility of replacing re-computation with analytical transformations.

Addressing the sensitivity of the method to parameters such as the agglomeration parameter, background mesh selection, and initial guess quality is another crucial aspect of our ongoing research. We aim to improve the method's robustness, ensuring stability even in the presence of minor parameter perturbations. Exploring alternative objective functions, including physics-based ones, and introducing terms into the objective function to penalize oscillatory level sets are avenues we are actively investigating.

\section*{Acknowledgements}
The work of Jakob Vandergrift is supported by the Graduate School CE within Computational Engineering at Technische Universität Darmstadt.
\section*{Conflict of interest} The corresponding author states on behalf of all authors, that there is no conflict of interest.
\bibliography{bibliography}
\bibliographystyle{ieeetr}

\end{document}